\documentstyle{amlts}
\begin{document}
\annalsline{157}{2003}
\received{September 24, 2001}
\startingpage{891}
\def\bye{\end{document}}
 \font\tenrm=cmr10
\def\ritem#1{\item[{\rm #1}]}
\input amssym.def
\input amssym.tex
\def\joinrel{\mathrel{\mkern-4mu}}
\def\relbar{\mathrel{\smash-}}
\def\lrar{\relbar\joinrel\relbar\joinrel\relbar\joinrel\relbar\joinrel\relbar\joinrel\relbar\joinrel\relbar\joinrel\relbar\joinrel\rightarrow}
\def\vlrar{\relbar\joinrel
\relbar\joinrel \relbar\joinrel \relbar\joinrel
\relbar\joinrel\relbar\joinrel
\relbar\joinrel\relbar\joinrel\relbar\joinrel\relbar\joinrel\relbar\joinrel\relbar\joinrel\relbar\joinrel\rightarrow}

\def\xrightarrow#1{\stackrel{#1}{\lrar}}
\def\vxrightarrow#1{\stackrel{#1}{\vlrar}}
\def\srightarrow#1{\stackrel{#1}{\longrightarrow}}
\catcode`\@=11
\font\twelvemsb=msbm10 scaled 1100
\font\tenmsb=msbm10
\font\ninemsb=msbm10 scaled 800
\newfam\msbfam
\textfont\msbfam=\twelvemsb  \scriptfont\msbfam=\ninemsb
  \scriptscriptfont\msbfam=\ninemsb
\def\msb@{\hexnumber@\msbfam}
\def\Bbb{\relax\ifmmode\let\next\Bbb@\else
 \def\next{\errmessage{Use \string\Bbb\space only in math
mode}}\fi\next}
\def\Bbb@#1{{\Bbb@@{#1}}}
\def\Bbb@@#1{\fam\msbfam#1}
\catcode`\@=12

 \catcode`\@=11
\font\twelveeuf=eufm10 scaled 1100
\font\teneuf=eufm10
\font\nineeuf=eufm7 scaled 1100
\newfam\euffam
\textfont\euffam=\twelveeuf  \scriptfont\euffam=\teneuf
  \scriptscriptfont\euffam=\nineeuf
\def\euf@{\hexnumber@\euffam}
\def\frak{\relax\ifmmode\let\next\frak@\else
 \def\next{\errmessage{Use \string\frak\space only in math
mode}}\fi\next}
\def\frak@#1{{\frak@@{#1}}}
\def\frak@@#1{\fam\euffam#1}
\catcode`\@=12
\newcommand{\bs}{\backslash}
\newcommand{\A}{{{\Bbb A}}}
\newcommand{\AAA}{{A}}
\newcommand{\AF}{{{\cal A}}}
\newcommand{\T}{{{\cal T}}}
\newcommand{\C}{{{\Bbb C}}}
\newcommand{\cmpct}{{{\cal C}}}
\newcommand{\Z}{{{\Bbb Z}}}
\newcommand{\Q}{{{\Bbb Q}}}
\newcommand{\K}{{{\bf K}}}
\newcommand{\U}{{{\Bbb U}}}
\newcommand{\SSS}{{{\cal S}}}
\newcommand{\BBB}{{{\cal B}}}
\newcommand{\CCC}{{{\cal C}}}
\newcommand{\DDD}{{{\cal D}}}
\newcommand{\WW}{{{\Bbb W}}}
\newcommand{\sprod}[2]{\left\langle#1,#2\right\rangle}
\newcommand{\inner}[2]{\left(#1,#2\right)}
\newcommand{\conj}[1]{\overline{#1}}
\newcommand{\Lang}{{\cal L}}
\newcommand{\R}{{{\Bbb R}}}
\newcommand{\abs}[1]{\left\vert#1\right\vert}
\newcommand{\norm}[1]{\left\Vert#1\right\Vert}
\newcommand{\projht}[1]{\left\vert#1\right\vert}
\newcommand{\M}{{{\frak M}}}
\newcommand{\Mabs}{{M}}
\newcommand{\m}{{m}}
\newcommand{\Outer}{\Sigma}
\newcommand{\sym}{{\mathop{\rm sym}}}
\newcommand{\expo}{{\mathop{\rm e}}}
\newcommand{\herm}{{{\frak B}}}
\newcommand{\least}{{{\cal L}}}
\newcommand{\jacquet}{{{\cal J}}}
\newcommand{\groth}{{{\frak R}}}
\newcommand{\sd}{{{\rm s.d.}}}
\newcommand{\sdu}{{{\rm s.d.u.}}}
\newcommand{\temp}{{{\rm temp}}}
\newcommand{\nontemp}{{{\rm n.t.}}}
\newcommand{\orth}{{{\rm non}\hbox{\scriptsize -}\SSS\hbox{\scriptsize
-}{\rm type}}}
\newcommand{\symp}{{{\rm pure}\hbox{\scriptsize -}\SSS\hbox{\scriptsize
-}{\rm type}}}
\newcommand{\allsymp}{{{\rm tot}\hbox{\scriptsize -}\SSS\hbox{\scriptsize
-}{\rm type}}}
\newcommand{\orthpairs}{{{\rm non}\hbox{\scriptsize -}\SSS\hbox{\scriptsize
-}{\rm pairs}}}
\newcommand{\red}{{red}}
\newcommand{\LQ}[1]{{{\cal LQ}(#1)}}             
\newcommand{\lqq}[1]{{{\cal SP}(#1)}}       
\newcommand{\dual}[1]{{\widetilde{#1}}}   
\newcommand{\opp}[1]{{{#1}^\flat}}       
\newcommand{\w}[1]{{{#1}^\sharp}}       
\newcommand{\hermdual}[1]{{{#1}^*}}       
\newcommand{\sdp}{{\rtimes}}
\newcommand{\nrmint}{{R}}
\newcommand{\intgln}{{{\cal M}}}
\newcommand{\nrmintt}{{{\cal R}}}
\newcommand{\Irr}{{\mathop{\rm Irr}}}
\newcommand{\resm}{{m_{-1}}}
\newcommand{\Inner}{{{\frak I}}}
\newcommand{\N}{{{\frak N}}}
\newcommand{\G}{{{\frak G}}}
\newcommand{\dummy}{\bullet}
\newcommand{\ident}{\iota}
\newcommand{\identt}{\kappa}
\newcommand{\V}{{{\Bbb V}}}
\newcommand{\conv}{{\star}}
\newcommand{\W}{{W}}
\newcommand{\triv}{{{\bf 1}}}
\newcommand{\id}{{{\bf 1}}}
\newcommand{\Ind}{{\mathop{\rm Ind}}}
\newcommand{\res}{{\mathop{\rm res}}}
\newcommand{\rad}{{\mathop{\rm rad}}}
\newcommand{\im}{{\mathop{\rm im}}}
\newcommand{\rank}{{\mathop{\rm rank}}}
\newcommand{\proj}{{\mathop{\rm proj}}}
\newcommand{\diag}{{\mathop{\rm diag}}}
\newcommand{\vol}{{\mathop{\rm vol}}}
\newcommand{\End}{{\mathop{\rm End}}}
\newcommand{\tr}{{\mathop{\rm tr}}\,}
\renewcommand{\Re}{{\mathop{\rm Re}}}
\renewcommand{\Im}{{\mathop{\rm Im}}}
\newcommand{\swrz}{{{\cal S}}}
\newcommand{\modulus}{{\delta}}
\newcommand{\eis}{{{\Bbb E}}}
\newcommand{\per}{{{\frak P}}}
\newcommand{\isom}{\simeq}
\newcommand{\form}{{{\cal F}}}
\newcommand{\reseis}{E_{-1}}
\newcommand{\resint}{{\M_{-1}}}
\newcommand{\rest}{{|}}
\newcommand{\Ht}{{H}}
\renewcommand{\a}{{{\frak a}}}
\newcommand{\eps}{{\varepsilon}}
\newcommand{\xprod}{{\mathsf{X}}} 
\newcommand{\Chi}{{\frak X}}
\newcommand{\sm}[4]{{\bigl(\begin{smallmatrix}{#1}&{#2}\\{#3}&{#4}
\end{smallmatrix}\bigr)}}
\newcommand{\Bessel}{{{\cal B}}}
\newcommand{\bil}{{{\frak B}}}

\font\emi= cmmi10 scaled 1700
\font\eightmi=cmmi10 
\font\titr=cmr10

\font\eighteenmsb=msbm10 scaled 1700
\font\eighteenrm=cmr10 scaled 1800
\title{On the nonnegativity of {\emi L}{\titr(}$\frac12,\pi${\titr)}  for
{\titr  SO}\lower2pt\hbox{$2${\eightmi n}$ +1$}}

 \def\titleheadline#1{\def\one{#1}\ifx\one\empty\else
\gdef\thetitle{{\frenchspacing%
\let\\ \relax
{#1}}}\fi}
\newif\ifshort
\def\shortname#1{\global\shorttrue\xdef
\theauthors{{\eightsc\uppercase{#1}}}}
\let\shorttitle\titleheadline
\shorttitle{\eightsc\uppercase{On the nonnegativity of} {\eightpoint \it
L}{\eightpoint
(}$\scriptscriptstyle\frac12,\pi${\eightpoint )} {\eightsc\uppercase{for}
{\eightpoint {\rm SO}}\lower2pt\hbox{2{\fivei
n}+1}}}
  \acknowledgements{First named author partially supported by NSF grant
DMS-0070611. Second named author
partially supported by NSF grant  DMS-9970342.}
 \twoauthors{Erez Lapid}{Stephen Rallis}
 \institutions{Einstein Institute of Mathematics, The Hebrew University of Jerusalem, Jerusalem  Israel\\
{\eightpoint {\it E-mail address\/}: erezla@math.huji.ac.il}\\
\vglue6pt
The Ohio
State University, Columbus,
 OH\\
{\eightpoint {\it E-mail address\/}: haar@math.ohio-state.edu}}

\centerline{\bf Abstract}
\vglue10pt
Let $\pi$ be a cuspidal generic representation of ${\rm SO}(2n+1,\A)$. We
prove that $L(\frac12,\pi)\ge0$.

\vglue-12pt

\section{Introduction}

Let $\pi$ be a cuspidal automorphic representation of ${\rm
GL}_n(\A)$ where $\A$ is the ring of ad\`eles of a number field
$F$. Suppose that $\pi$ is self-dual. Then the ``standard''
$L$-function (\cite{GJ72}) $L(s,\pi)$ is real for $s\in\R$ and
positive for $s>1$. Assuming GRH we have $L(s,\pi)>0$ for
$\frac12<s\le1$, except for the case where $n=1$ and $\pi$ is the
trivial character. It would follow that $L(\frac12,\pi)\ge0$.
However, the latter is not known even in the case of quadratic
Dirichlet characters. In general, if $\pi$ is self-dual then $\pi$
is either \emph{symplectic} or \emph{orthogonal}, i.e. exactly one
of the (partial) $L$-functions $L^S(s,\pi,\wedge^2)$,
$L^S(s,\pi,\sym^2)$ has a pole at $s=1$. In the first case $n$ is
even and the central character of $\pi$ is trivial (\cite{JS90a}).
In the language of the Tannakian formalism of Langlands
(\cite{Lan79}), any cuspidal representation $\pi$ of ${\rm
GL}_n(\A)$ corresponds to an irreducible $n$-dimensional
representation $\varphi$ of a conjectural group $\Lang_F$ whose
derived group is compact. Then $\pi$ is self-dual if and only if
$\varphi$ is self-dual, and the classification into symplectic and
orthogonal is compatible with (and suggested by) the one for
finite dimensional representations of a compact group. Our goal in
this paper is to show

\specialnumber{1}\proclaim{Theorem} \label{Main}
Let $\pi$ be a symplectic cuspidal representation of ${\rm GL}_n(\A)$.
Then $L(\frac12,\pi)\ge0$.
\endproclaim

We note that the same will be true for the partial $L$-function.
The value $L(\frac12,\pi)$ appears in many arithmetic, analytic
and geometric contexts -- among them, the Shimura correspondence
(\cite{Wal81}), or more generally -- the theta
correspondence (\cite{Ral87}), the Birch-Swinnerton-Dyer
conjecture, the Gross-Prasad conjecture (\cite{GP94}),
certain period integrals, and the relative trace formula
(\cite{JC01}, \cite{BM}). In all the above cases, the
$L$-functions are of symplectic type. Moreover, all motivic
$L$-functions which have the center of symmetry as a critical
point in the sense of Deligne are necessarily of symplectic type.
In the case $n=2$, $\pi$ is symplectic exactly when the central
character of $\pi$ is trivial. The above-mentioned interpretations
of $L(\frac12,\pi)$ were used to prove Theorem \ref{Main} in that
case (\cite{KZ81}, \cite{KS93}, using the Shimura
correspondence in special cases, and \cite{Guo96}, using a
variant of Jacquet's relative trace formula, in general). The
nonnegativity of $L(\frac12,\pi)$ in the ${\rm GL}_2$ case already has
striking applications, for example to sub-convexity estimates for
various $L$-functions (\cite{CI00}, \cite{Ivi01}). We
expect that the higher rank case will turn out to be useful as
well. The nonnegativity of $L(\frac12,\chi)$ for quadratic
Dirichlet characters would have far-reaching implications to Gauss
class number problem. Unfortunately, our method is not applicable
to that case.

The Tannakian formalism suggests that the symplectic and
orthogonal automorphic representations of ${\rm GL}_n(\A)$ are
functorial images from classical groups. In fact, it is known that
every symplectic cuspidal automorphic representation $\pi$ of
${\rm GL}_{2n}(\A)$ is a functorial image of a cuspidal generic
representation of ${\rm SO}(2n+1,\A)$. Conversely, to every cuspidal
generic representation of ${\rm SO}(2n+1,\A)$ corresponds an automorphic
representation of ${\rm GL}_{2n}(\A)$ which is parabolically induced
from cuspidal symplectic representations
(\cite{GRS01},\break \cite{CKP-SS01}). As a consequence:

\specialnumber{2}\proclaim{Theorem} \label{Mainref}
Let $\sigma$ be a cuspidal generic representation of
${\rm SO}(2n+1,\A)$. Then $L^S\left(\frac12,\sigma\right)\ge0$.
\endproclaim

The $L$-function is the one pertaining to the imbedding of ${\rm
Sp}(n,\C)$, the\break $L$-group of ${\rm SO}(2n+1)$, in ${\rm
GL}(2n,\C)$. By the work of Jiang-Soudry (\cite{JS}) Theorem
\ref{Mainref} applies equally well to the completed $L$-function
as defined by Shahidi in \cite{Sha81}.

We emphasize however that our proof of Theorem \ref{Main} is
independent of the functorial lifting above. In fact, it turns
out, somewhat surprisingly, that Theorem \ref{Main} is a simple
consequence of the theory of Eisenstein series on classical
groups. Consider the symplectic group ${\rm Sp}_n$ and the
Eisenstein series $E(g,\varphi,s)$ induced from $\pi$ viewed as a
representation on the Siegel parabolic subgroup. If $\pi$ is
symplectic then for $E(g,\varphi,s)$ to have a pole at $s=\frac12$
it is necessary and sufficient that $L(\frac12,\pi)\ne0$, in which
case the pole is simple. In particular, in this case
$\eps(\frac12,\pi)=1$ by the functional equation. We refer the
reader to the body of the paper for any unexplained notation. Let
$\reseis(\cdot,\varphi)$ be the residue of $E(\cdot,\varphi,s)$ at
$s=\frac12$. It is a square-integrable automorphic form on ${\rm
Sp}_n$. A consequence of the spectral theory is that the inner
product of two such residues is given by the residue $\resint$ of
the intertwining operator at $s=\frac12$. Thus, $\resint$ is a
positive semi-definite operator. First assume that the local
components of $\pi$ are unramified at every place including the
archimedean ones. Then by a well-known formula of Langlands
(\cite{Lan71}), the intertwining operator $\M(s)$ satisfies
$$
\M(s)v_0=L(s,\pi)/L(s+1,\pi)\cdot
L(2s,\pi,\wedge^2)/L(2s+1,\pi,\wedge^2)\cdot v_0
$$
for the unramified vector $v_0$. Therefore
$$
\resint v_0=\frac12\cdot
L\left(\frac12,\pi\right)/L\left(\frac32,\pi\right)\cdot\res_{s=1}
L(s,\pi,\wedge^2)/L(2,\pi,\wedge^2)\cdot v_0.
$$
Since $L(s,\pi)$ is positive for $s>1$ and $L(s,\pi,\wedge^2)$ is
real and nonzero for $s>1$ we obtain Theorem \ref{Main} in this
case. In order to generalize this argument and avoid any local
assumptions on $\pi$ we have, as usual, to make some local
analysis. For that, we use Shahidi's normalization of the
intertwining operators (\cite{Sha90b}) which is applicable since
$\pi$ is generic. Let
$\nrmint(\pi,s)=\nrmint(s)=\otimes_v\nrmint_v(s):I(\pi,s)\rightarrow
I(\pi,-s)$ be the normalized intertwining operator. Here we take
into account a canonical identification of $\pi$ with its
contragredient and suppress the dependence of $\nrmint_v(s)$ on a
choice of an additive character. Then $\M(s)=\m(s)\cdot\nrmint(s)$
where
$$
\m(s)=\frac{L(s,\pi)}{\eps(s,\pi)L(s+1,\pi)}\cdot
\frac{L(2s,\pi,\wedge^2)}{\eps(2s,\pi,\wedge^2)L(2s+1,\pi,\wedge^2)}.
$$
Hence, $\resint=\resm\cdot\nrmint\left(\frac12\right)$, where
$\resm$ is the residue of $\m(s)$ at $s=\frac12$, and the operator
$\nrmint(\frac12)$ is semi-definite with the same sign as $\resm$.
On the other hand, the argument of Keys-Shahidi (\cite{KS88})
shows that the Hermitian involution $\nrmint(\pi_v,0)$ has a
nontrivial $+1$ eigenspace. The main step (Lemma \ref{claim2},
proved in \S\ref{locanal}) is to show that
$\nrmint(\pi_v,\frac12)$ is \emph{positive} semi-definite by
``deforming'' it to $\nrmint(\pi_v,0)$. This will imply that
$\resm>0$, i.e.
$$
\frac{L\left(\frac12,\pi\right)}{L\left(\frac32,\pi\right)}\cdot
\frac{\res_{s=1}L(s,\pi,\wedge^2)}{\eps(1,\pi,\wedge^2)L(2,\pi,\wedge^2)}>0.
$$
Similarly, working with the group ${\rm SO}(2n)$ we obtain
$$
\frac{\res_{s=1}L(s,\pi,\wedge^2)}{\eps(1,\pi,\wedge^2)L(2,\pi,\wedge^2)}>0
$$
if $\pi$ is symplectic. Altogether this implies Theorem \ref{Main}
(see \S\ref{setup}). We may work with the group ${\rm SO}(2n+1)$
as well. Using the relation $\eps(\frac12,\pi\otimes\dual{\pi})=1$
(\cite{BH99}) we will obtain the following:

\specialnumber{3}\proclaim{Theorem} \label{Main2} Let $\pi$ be a
self\/{\rm -}\/dual cuspidal representation of ${\rm GL}_n(\A)$.
Then
$\eps(\frac12,\pi,\wedge^2)=\eps(\frac12,\pi,\sym^2)=1$.
\endproclaim

This is compatible with the Tannakian formalism. In general one
expects that $\eps(\frac12,\pi,\rho)=1$ if the representation
$\rho\circ\varphi$ is orthogonal (\cite{PR99}). This is
inspired by results of Fr{\"o}hlich-Queyrut, Deligne and Saito
about epsilon factors of orthogonal Galois representations and
motives (\cite{FQ73}, \cite{Del76}, \cite{Sai95}).

The analysis of Section~\ref{locanal}, the technical core of this
article, relies on detailed information about the reducibility of
induced representations of classical groups. This was studied
extensively by Goldberg, Jantzen, Muic, Shahidi, Tadic, and others
(see \cite{Gol94}, \cite{Jan96}, \cite{Mui01},
\cite{Sha92}, \cite{Tad98}).

\demo{Note added in proof} Since the time of writing this paper
Theorem \ref{Main} was generalized by the first-named author to
tensor product $L$-functions of symplectic type (\cite{Lap03}).
Similarly, other root numbers of orthogonal type have shown to be
$1$ (\cite{Lap02}).

The authors would like to express their gratitude to the Institute
for Advanced Study for the hospitality during the first half of
2001. We would also like to thank Professors Herv\'e Jacquet and
Freydoon Shahidi for useful discussions.

\section{The setup} \label{setup}

Let $F$ be a number field, $\A=\A_F$ its ad\`eles ring and let
$\pi$ be a cuspidal automorphic representation of ${\rm
GL}_n(\A)$. We say that $\pi$ is symplectic (resp. orthogonal) if
$L^S(s,\pi,\wedge^2)$ (resp. $L^S(s,\pi,\sym^2)$) has a pole at
$s=1$. If $\pi$ is symplectic or orthogonal then $\pi$ is
self-dual. Conversely, if $\pi$ is self-dual then $\pi$ is either
symplectic or orthogonal but not both. Moreover, if $\pi$ is
symplectic then $n$ is even and the central character of $\pi$ is
trivial (\cite{JS90a}). Our goal is to prove Theorems \ref{Main}
and \ref{Main2}. In this section we will reduce them to a few
local statements, namely Lemmas \ref{local}--\ref{rankin} below
which will be proved in the next section. They all have some
overlap with known results in the literature. We first fix some
notation. By our convention, if $X$ is an algebraic group over $F$
we denote the $F$-points of $X$ by $X$ as well. Let $J_n$ be
$n\times n$ matrix with ones on the nonprincipal diagonal and
zeros otherwise. Let $G$ be either the split orthogonal group
${\rm SO}(2n+1)$ with respect to the symmetric form defined by
$$
\left(\begin{array}{ccc}&&J_n\\&1&\\J_n&&\end{array}\right)
$$
or the symplectic group ${\rm Sp}_n$ with respect to the
skew-symmetric form defined by the matrix
$$
\left(\begin{array}{cc}0&J_n\\-J_n&0\end{array}\right)
$$
or the split orthogonal group ${\rm SO}(2n)$ with respect to the
symmetric form defined by
$\left(\begin{array}{cc}0&J_n\\J_n&0\end{array}\right)$. Then $G$
acts by right multiplication on the space $\V$ of row vectors of
size $2n$ or $2n+1$. Let $P=M\cdot U$ be the Siegel parabolic
subgroup of $G$ with its standard Levi decomposition. It is the
stabilizer of the maximal isotropic space $\U$ defined by the
vanishing of all but the last $n$ coordinates. We identify $M$
with ${\rm GL}(\V/\U^\perp)\isom {\rm GL}_n$ where $\U^\perp$ is
the perpendicular of $\U$ in $\V$ with respect to the form
defining $G$. We denote by $\nu:M(\A)\rightarrow\R_+$ the absolute
value of the determinant in that identification. Let $\K$ be the
standard maximal compact subgroup of $G(\A)$. We extend $\nu$ to a
left-$U(\A)$ right-$\K$-invariant function on $G(\A)$ using the
Iwasawa decomposition. Let $\modulus_P$ be the modulus function of
$P(\A)$. It is given by $\delta_P=\nu^n,\nu^{n+1}$ or $\nu^{n-1}$
according to whether $G={\rm SO}(2n+1)$, ${\rm Sp}_n$ or ${\rm
SO}(2n)$. Let $\pi$ be a cuspidal representation of ${\rm
GL}_n(\A)$ and $\AF(U(\A)M\bs G(\A))_{\pi,s}$ be the space of
automorphic forms $\varphi$ on $U(\A)M\bs G(\A)$ such that the
function $m\rightarrow\nu^{-s}(m)\modulus_P(m)^{-1/2}\varphi(mk)$
belongs to the space of $\pi$ for any $k\in\K$. By multiplicity-one for ${\rm GL}_n$, $\AF(U(\A)M\bs G(\A))_{\pi,s}$
depends only on the equivalence class of $\pi$ and not on its automorphic
realization. By choosing an automorphic realization for $\pi$
(unique up to a scalar), we may identify $\AF(U(\A)M\bs
G(\A))_{\pi,s}$ with (the $\K$-finite vectors in) the induced
space $I(\pi,s)$. The Eisenstein series
$$
E(g,\varphi,s)=\sum_{\gamma\in P\bs G}\varphi(\gamma
g)\nu^s(\gamma g)
$$
converges when $\Re(s)$ is sufficiently large and admits a
meromorphic continuation. Whenever it is regular it defines an
intertwining map $\AF(U(\A)M\bs G(\A))_{\pi,s}\break\rightarrow\AF(G\bs
G(\A))$. It is known that the only possible singularity of
$E(g,\varphi,s)$ for $\Re(s)\ge0$ is a simple pole at $s=\frac12$
(except when $\pi$ is the trivial character and $G={\rm Sp}_1$, where
there is a pole at $s=1$).

In the case $G={\rm SO}(2n)$ let $\Outer$ be the outer automorphism
obtained by conjugation by the element
$$
\left(\begin{array}{cccc}1_{n-1}&&&\\&0&1&\\&1&0&\\&&&1_{n-1}\end{array}\right)
$$
of ${\rm O}(2n)\setminus {\rm SO}(2n)$.
For the other groups let $\Outer=\id$. In all cases we set
$\theta=\Outer^n$. Then $\theta$ induces the principal involution
on the root data of $G$. Note that $\{P,\theta(P)\}$ is the set of
standard parabolic subgroups of $G$ which are associate to $P$.
Fix $w\in G\setminus M$ such that $wMw^{-1}=\theta(M)$; it is
uniquely determined up to right multiplication by $M$. Let
$\w{}:M\rightarrow\theta(M)$ be defined by $\w{m}=wmw^{-1}$.
Denote by $w\pi$ the cuspidal automorphic representation of
$\theta(M)(\A)$ on $\{\w{\varphi}:\varphi\in V_\pi\}$ where
$\w{\varphi}(\w{m})=\varphi(m)$. The ``automorphic'' intertwining
operator
$$
\M(s)=\M(\pi,s):\AF(U(\A)M\bs
G(\A))_{\pi,s}\rightarrow\AF(\theta(U)(\A)\theta(M)\bs
G(\A))_{w\pi,-s}
$$
is defined by
$$
[\M(s)\varphi](g)=\int_{\theta(U)(\A)}\varphi(w^{-1}ug)\nu^s(w^{-1}ug)\
du.
$$
Let $\reseis(\dummy,\varphi)$ be the residue of $E(g,\varphi,s)$
at $s=\frac12$. It is zero unless $w\pi=\pi$, and in particular,
$\theta(M)=M$, i.e. $\theta=\id$. The latter means that $P$ is
conjugate to its opposite. We say that $\pi$ is of $G$-type if
$\reseis\not\equiv0$, or what amounts to the same, that
$\resint\not\equiv0$ where $\resint$ is the residue of $\M(s)$ at
$\frac12$. In this case $\reseis$ defines an intertwining map
$\AF(U(\A)M\bs G(\A))_{\pi,\frac12}\rightarrow\AF(G\bs G(\A))$.
The inner product formula for two residues of Eisenstein series is
given by
\begin{eqnarray} \label{inner}
&&\hskip-36pt\int_{G\bs
G(\A)}\reseis(g,\varphi_1)\conj{\reseis(g,\varphi_2)}\
dg\\
&&\hskip.5in =\int_{\K}\int_{M\bs
M(\A)^1}\resint\varphi_1(mk)\conj{\varphi_2(mk)}\ dm\ dk\nonumber
\end{eqnarray}
up to a positive constant depending on normalization of Haar
measures. This follows for example by taking residues in the
Maass-Selberg relations for inner product of truncated Eisenstein series (%
cf.\ \cite[\S4]{Art80}). Alternatively, this is a consequence
of spectral theory (\cite{MW95}).

We let $\w{\pi}$ be the representation of $\theta(M)(\A)$ on
$V_\pi$ defined by $\w{\pi}(\w{m})v=\pi(m)v$. We may identify
$\w{\pi}$ with $w\pi$ by the map $\varphi\mapsto\w{\varphi}$. Let
$\Mabs(s)=\Mabs(\pi,s):I(\pi,s)\rightarrow I(\w{\pi},-s)$ be the
``abstract'' intertwining operator given by
$$
\Mabs(s)\varphi(g)=\int_{\theta(U)(\A)}\varphi(w^{-1}ug)\nu^s(w^{-1}ug)\
du.
$$
Under the isomorphisms
\begin{eqnarray*}
\AF(U(\A)M\bs G(\A))_{\pi,s}&\isom&
I(\pi,s) \hbox{ and}\\
 \AF(\theta(U)(\A)\theta(M)\bs G(\A))_{w\pi,-s}&\isom&
I(\w{\pi},-s),
\end{eqnarray*} $\M(s)$ becomes $\Mabs(s)$.

Let $\opp{}:M\rightarrow M$ be the map defined by
$\opp{m}=\theta(\w{m})$. We will choose the representative $w$ as in
\cite{Sha90b} so that when $M$ is identified with ${\rm GL}_n$,
$\opp{}$ becomes the involution $x\mapsto w_n^{-1}\,^tx^{-1}w_n$ where
$$
(w_n)_{ij}=\left\{ \begin{array}{ll}(-1)^i&\hbox{if }i+j=n+1\\0
&\hbox{otherwise.}\end{array}\right.
$$
In particular $\opp{}$ does not depend on $G$. A direct computation
shows that
\begin{equation} \label{w2}
w^2\in M\textrm{ corresponds to the central element }(-1)^n
\textrm{(resp. $(-1)^{n+1}$) of } GL_n
\end{equation}
if $G$ is symplectic (resp. orthogonal).
We define
$\opp{\varphi}$ and $\opp{\pi}$ as before. Since $\pi$ is
irreducible we have (\cite{GK75})
\begin{equation} \label{contra}
\opp{\pi}\hbox{ is equivalent to the contragredient }
\dual{\pi}\hbox{ of }\pi.
\end{equation}
Thus, for $\pi$ to be of $G$-type it is necessary that
$\theta=\id$ and that $\pi$ be self-dual. If $\pi$ is self-dual we
define the intertwining operator
$\ident=\ident_\pi:\opp{\pi}\rightarrow\pi$ by
$\ident(\varphi)=\opp{\varphi}$. It is well-defined by
multiplicity-one and does not depend on the automorphic
realization of $\pi$. We write $\ident(s)=\ident(\pi,s)$ for the
induced map $I(\opp{\pi},s)\rightarrow I(\pi,s)$ given by
$\left[\ident(s)(f)\right](g)=\ident(f(g))$. Note that when
$\theta=\id$, $\ident(s)$ is the map $I(\opp{\pi},s)\rightarrow
I(\pi,s)$ induced from the ``physical'' equality of the two spaces
$\AF(U(\A)M\bs G(\A))_{w\pi,s}$ and $\AF(U(\A)M\bs
G(\A))_{\pi,s}$. Assume that $\pi$ is self-dual and that
$\theta=\id$. Then as a map from $I(\pi,s)$ to $I(\pi,-s)$ the
intertwining operator $\M(s)$ becomes $\ident(-s)\circ\Mabs(s)$.
Let $(\cdot,\cdot)_\pi$ be the invariant positive-definite
Hermitian form on $\pi$ obtained through its automorphic
realization. This gives rise to the invariant sesqui-linear form
$(\cdot,\cdot)=(\cdot,\cdot)_s:I(\pi,-s)\times
I(\pi,\conj{s})\rightarrow\C$ given by
$$
(\varphi_1,\varphi_2)=\int_{\K}(\varphi_1(k),\varphi_2(k))_\pi\
dk.
$$
Thus, the right-hand side of (\ref{inner}), viewed as a
positive-definite invariant Hermitian form on $I(\pi,\frac12)$, is
$(\ident(-\frac12)\circ\Mabs_{-1}\varphi_1,\varphi_2)_{\frac12}$.

In the local case we can define $\w{\pi_v}$, $\opp{\pi_v}$ and the
local intertwining operators
$$
\Mabs_v(s):I(\pi_v,s)\rightarrow I(\w{\pi_v},-s)
$$
in the same way. Fix a nontrivial character $\psi=\otimes_v\psi_v$ of
$F\bs\A_F$. For any $v$ choose a Whittaker model for $\pi_v$ with respect
to the $\opp{}$-stable character
$$
\left(\begin{array}{cccc}1&x_1&*&*\\&1&\ddots&*\\&&1&x_{n-1}\\&&&1\end{array}\right)
\mapsto\psi_v(x_1+\dots+x_{n-1}).
$$
If $\pi_v$ is self-dual then we define the intertwining map
$\ident_v=\ident^{\psi_v}_{\pi_v}:\opp{\pi _v}\rightarrow\pi_v$ by
$$
\left[\ident_v(W)\right](g)=W(\opp{g})
$$
in the Whittaker model with respect to $\psi_v$.
By uniqueness of the  Whittaker model
$\ident_v$ is well-defined and does not depend on  choice of
the Whittaker model. If we change $\psi_v$ to $\psi_v(a\cdot)$
for $a\in F_v^*$ then $\ident_v$ is multiplied by the sign
$\omega_{\pi_v}^{n-1}(a)$.
If $\pi_v$ and $\psi_v$ are unramified then $\ident_v(u)=u$
for an unramified
vector $u$ since the unramified Whittaker vector is nonzero at
the identity by the Casselman-Shalika formula.

Suppose that $\pi=\otimes_v\pi_v$ is an automorphic self-dual
cuspidal representation of ${\rm GL}_n(\A)$ where the restricted
tensor product is taken with respect to a choice of unramified
vectors $e_v$ almost everywhere. We choose invariant positive
definite Hermitian forms $(\cdot,\cdot)_{\pi_v}$ on $\pi_v$ for
all $v$ so that $(e_v,e_v)_{\pi_v}=1$ almost everywhere. This
gives rise to sesqui-linear forms
$(\cdot,\cdot)_{v,s}:I(\pi_v,-s)\times
I(\pi_v,\conj{s})\rightarrow\C$ as above. We have
$(\cdot,\cdot)_\pi=c\otimes(\cdot,\cdot)_{\pi_v}$ and
$(\cdot,\cdot)_s=c\otimes(\cdot,\cdot)_{v,s}$ in the obvious
sense, for some \emph{positive} scalar $c$, and
$\ident_\pi=\otimes_v\ident_{\pi_v}$. \vglue4pt

At this point it is useful to normalize $\Mabs_v(s)$ by the
normalization factors $\m_v^{\psi_v}(\pi_v,s)=\m_v(s)$ defined by Shahidi
in \cite{Sha90b}. The latter are given by
$$
\m_v(s)=\left\{
\begin{array}{ll}\frac{L(2s,\pi_v,\sym^2)}{\eps(2s,\pi_v,\sym^2,\psi_v^{-1})
L(2s+1,\pi_v,\sym^2)}&G={\rm SO}(2n+1),\\[6pt]
\frac{L(s,\pi_v)}{\eps(s,\pi_v,\psi_v^{-1})L(s+1,\pi_v)}\frac{L(2s,\pi_v,\wedge^2)}{
\eps(2s,\pi_v,\wedge^2,\psi_v^{-1})L(2s+1,\pi_v,\wedge^2)}&G={\rm Sp}_n,\\[6pt]
\frac{L(2s,\pi_v,\wedge^2)}{
\eps(2s,\pi_v,\wedge^2,\psi_v^{-1})L(2s+1,\pi_v,\wedge^2)}&G={\rm SO}(2n),\end{array}\right.
$$
where $L(s,\pi_v)$, $L(s,\pi_v,\wedge^2)$, $L(s,\pi_v,\sym^2)$ are
the local $L$-functions pertaining to the standard, symmetric
square and exterior square representations of $GL_n(\C)$
respectively, and similarly for the epsilon factors. We write
$\Mabs_v(\pi_v,s)=\m_v^{\psi_v}(\pi_v,s)\nrmint_v^{\psi_v}(\pi_v,s)$
where $\nrmint_v(s)=\nrmint_v^{\psi_v}(\pi_v,s)$ are the
normalized intertwining operators. Note that by changing $\psi_v$
to $\psi_v(a\cdot)$ the scalar $\m_v(s)$ is multiplied by
$(\omega_{\pi_v}(a)\abs{a}^{n(s-\frac12)})^k$ where $k=n+1, n$, or
$n-1$ according to whether $G={\rm SO}(2n+1), {\rm Sp}_n$ or ${\rm
SO}(2n)$.

The following lemma will be proved in the next section, together with the
other lemmas below.

\specialnumber{1}\proclaim{Lemma} \label{local} For
all $v${\rm ,} $\nrmint_v(s)${\rm ,} $\Mabs_v(s)${\rm ,}
$L_v(2s,\pi_v,\sym^2)${\rm ,} $L_v(2s,\pi_v,\wedge^2)${\rm
,}\break $L_v(s,\pi_v)$ and $\m_v(s)$ are holomorphic and nonzero
for $\Re(s)\ge\frac12$.
\endproclaim

In fact, the holomorphy and nonvanishing of $\nrmint_v(s)$ for
$\Re(s)\ge\frac12$ is proved more generally in a recent paper of
Kim (\cite{Kim02}).

Let $\m(s)=\m(\pi,s)=\prod_v\m_v^{\psi_v}(\pi_v,s)$ and
$\nrmint(s)=\otimes_v\nrmint_v(s)$ so that
$\Mabs(s)=\m(s)\nrmint(s)$. If $G={\rm SO}(2n+1)$ then
$$
\m(s)=\frac{L(2s,\pi,\sym^2)}{\eps(2s,\pi,\sym^2)L(2s+1,\pi,\sym^2)}
=\frac{L(1-2s,\pi,\sym^2)}{L(1+2s,\pi,\sym^2)}.
$$
If $G={\rm Sp}_n$ then
\begin{eqnarray*}
\m(s)&=&\frac{L(s,\pi)}{\eps(s,\pi)L(s+1,\pi)}\frac{L(2s,\pi,\wedge^2)}{
\eps(2s,\pi,\wedge^2)L(2s+1,\pi,\wedge^2)}\\&=&\frac{L(1-s,\pi)}{L(1+s,\pi)}
\frac{L(1-2s,\pi,\wedge^2)}{L(1+2s,\pi,\wedge^2)}.
\end{eqnarray*}
If $G={\rm SO}(2n)$,
$$
\m(s)=\frac{L(2s,\pi,\wedge^2)}{
\eps(2s,\pi,\wedge^2)L(2s+1,\pi,\wedge^2)}=
\frac{L(1-2s,\pi,\wedge^2)}{L(1+2s,\pi,\wedge^2)}.
$$
In particular, the residue $\resm$ at $s=\frac12$ is equal to
$\frac12$ times
$$
\left\{\begin{array}{ll}\frac{\res_{s=1}L(s,\pi,\sym^2)}
{\eps(1,\pi,\sym^2)L(2,\pi,\sym^2)}&G={\rm SO}(2n+1)\\[6pt]
\frac{L(\frac12,\pi)}{\eps(\frac12,\pi)L(\frac32,\pi)}
\frac{\res_{s=1}L(s,\pi,\wedge^2)}
{\eps(1,\pi,\wedge^2)L(2,\pi,\wedge^2)}&G={\rm Sp}_n\\[6pt]
\frac{\res_{s=1}L(s,\pi,\wedge^2)}
{\eps(1,\pi,\wedge^2)L(2,\pi,\wedge^2)}&G={\rm
SO}(2n).\end{array}\right.
$$
By Lemma \ref{local}, $\pi$ is of $G$-type if and only if $m(s)$
has a pole (necessarily simple) at $s=\frac12$. Thus, $\pi$ is of
${\rm Sp}_n$ type if and only if $\pi$ is symplectic and
$L(\frac12,\pi)\ne0$; $\pi$ is of ${\rm SO}(2n+1)$ type if and
only if $\pi$ is orthogonal; $\pi$ is of ${\rm SO}(2n)$ type if
and only if $\pi$ is symplectic. Suppose that $\pi$ is of
$G$-type. Let $\herm(s)=\herm(\pi,s)$ be the operator
$\ident(-s)\circ\nrmint(s):I(\pi,s)\rightarrow I(\pi,-s)$ for
$s\in\R$ and let $\Inner(\pi,s)$ be the form on $I(\pi,s)$ defined
by $(\herm(s)\varphi,\varphi)$. Since
$\resint=\resm\cdot\herm\left(\frac12\right)$, it follows from
(\ref{inner}) that $\Inner(\pi,\frac12)$ is semi-definite with the
same sign as $\resm$. We will show that
\begin{equation} \label{show}
\Inner(\pi,\frac12)\hbox{ is positive semi-definite}
\end{equation}
and thus
\begin{equation} \label{resm}
\resm>0.
\end{equation}
\vglue12pt \demo{{\rm 2.1.} Proof of Theorem {\rm \ref{Main}}} We
will use (\ref{resm}) for the groups ${\rm Sp}_n$ and ${\rm
SO}(2n)$. Together, this implies that if $\pi$ is symplectic and
$L(\frac12,\pi)\ne0$ then
$\frac{L(\frac12,\pi)}{\eps(\frac12,\pi)L(\frac32,\pi)}\break >0$.
By the functional equation and the fact that $L(\frac12,\pi)\ne0$
we must have $\eps(\frac12,\pi)=1$. On the other hand $L(s,\pi)$
is a convergent Euler product for $s>1$, all factors of which are
real and positive. Indeed, $L(s,\pi_v)=\overline{L(\bar s,\pi_v)}$
since $\pi_v$ is equivalent to its Hermitian dual. In the
nonarchimedean case, $L(s,\pi_v)\rightarrow 1$ as
$s\rightarrow+\infty$ ($s$ real). In the archimedean case
$L(s,\pi_v)=\prod_{i=1}^n\Gamma_{\R}(s-s_i)$ for some $s_i\in\C$
where $\Gamma_{\R}(s)=\pi^{-s/2}\Gamma(s/2)$. We have $\sum\Im
s_i=0$ since $\pi_v=\overline{\pi}_v$. It is easily deduced from
Stirling's formula that $L(s,\pi_v)\rightarrow+\infty$ as
$s\rightarrow+\infty$. In both cases $L(s,\pi_v)$ is holomorphic
and nonzero for $s\ge\frac12$. The claim follows. Hence
$L(\frac32,\pi)>0$, and therefore, $L(\frac12,\pi)>0$.

It remains to prove (\ref{show}). The operator $\herm(\pi,s)$ and
the form $\Inner(\pi,s)$ admit a local analogue and we have
$\herm(\pi,s)=\otimes_v\herm^{\psi_v}(\pi_v,s)$ and
$\Inner(\pi,s)=c\otimes_v\Inner^{\psi_v}(\pi_v,s)$.

We will prove the following purely local Lemmas. Recall the
assumption that $\theta=\id$.
\enddemo

\specialnumber{2}\proclaim{Lemma} \label{claim} Let $\pi_v$ be a
generic irreducible unitary self\/{\rm -}\/dual representation of
${\rm GL}_n$ over a local field of characteristic $0$. Then
$\herm^{\psi_v}(\pi_v,s)$ is Hermitian for $s\in\R$ and
holomorphic near $s=0$. Moreover{\rm ,} $\herm^{\psi_v}(\pi_v,0)$ is an
involution with a nontrivial $+1$\/{\rm -}\/eigenspace.
\endproclaim

\specialnumber{3}\proclaim{Lemma} \label{claim2} Under the same
assumptions{\rm ,} suppose further that $\Inner^{\psi_v}(\pi_v,\frac12)$ is
semi\/{\rm -}\/definite. Then $\Inner^{\psi_v}(\pi_v,0)$ is definite with
the same sign as $\Inner^{\psi_v}(\pi_v,\frac12)$. Hence{\rm ,} by Lemma
{\rm \ref{claim},} $\herm^{\psi_v}(\pi_v,0)=\id$ and
$\Inner^{\psi_v}(\pi_v,\frac12)$ is positive semi\/{\rm -}\/definite.
\endproclaim

These two lemmas, together with the fact that
$\Inner(\pi,\frac12)$ is semi-definite, imply (\ref{show}), even
locally.

We remark that in the case where $G$ is an orthogonal group then
up to a positive scalar $\herm^{\psi_v}(\pi_v,s)$ is independent
of $\psi_v$. This is no longer true in the $Sp_n$ case if the
central character of $\pi_v$ is nontrivial. In that case, Lemma
\ref{claim} actually implies the well-known fact that $I(\pi_v,0)$
is reducible.

Note also that the very last (and most important) conclusion of
Lemma \ref{claim2} is trivial in the unramified case. Finally, let
us mention that a property related (and ultimately, equivalent) to
the conclusion of Lemma \ref{claim2} for the local components of a
symplectic cuspidal representation was proved by Jiang-Soudry
using the descent construction (\cite{JS}). We will not use their
result.

\demo{{\rm 2.2.} Proof of Theorem \ref{Main2}} We first observe
that $L(s,\pi,\sym^2)$ and\break $L(s,\pi,\wedge^2)$ are
holomorphic and nonzero for $\Re(s)>1$. Indeed, the partial\break
$L$-functions $L^S(s,\pi,\sym^2)$, $L^S(s,\pi,\wedge^2)$ are
holomorphic for $\Re(s)>1$ (\cite{JS90a}, \cite{BG92}) and their
product is $L^S(s,\pi\otimes\pi)$, which is nonzero for
$\Re(s)>1$, since the Euler product converges absolutely
(\cite{JS81}). The statement now follows from Lemma~\ref{local}.

Suppose that $\pi$ is orthogonal. Applying (\ref{resm}) to the
group ${\rm SO}(2n+1)$ we obtain $\frac{\res_{s=1}L(s,\pi,\sym^2)}
{\eps(1,\pi,\sym^2)L(2,\pi,\sym^2)}>0$. Since $L(s,\pi,\sym^2)$ is
real and nonzero for  $s>~1$ we obtain
$\frac{\res_{s=1}L(s,\pi,\sym^2)} {L(2,\pi,\sym^2)}>0$. Hence
$\eps(1,\pi,\sym^2)>0$. Since $\eps(s,\pi,\sym^2)$ is nonzero and
real for $s\in\R$ we get $\eps(\frac12,\pi,\sym^2)>0$. On the
other hand, $\eps(\frac12,\pi,\sym^2)=\pm1$ by the functional
equation  and hence, $\eps(\frac12,\pi,\sym^2)=1$. Similarly, if
$\pi$ is symplectic then using the group $G={\rm SO}(2n)$ and the
same argument we obtain $\eps(\frac12,\pi,\wedge^2)=1$. Since any
self-dual cuspidal representation $\pi$ is either symplectic or
orthogonal, the above argument shows that either
$\eps(\frac12,\pi,\wedge^2)=1$ or $\eps(\frac12,\pi,\sym^2)=1$. On
the other hand for any $\pi$ (self-dual or not)
\begin{equation} \label{epseps}
\eps(s,\pi\otimes\pi)=\eps(s,\pi,\wedge^2)\eps(s,\pi,\sym^2).
\end{equation}
Indeed, this follows from the corresponding equality of
$L$-functions, which is in fact true locally. In the archimedean
case this follows from the compatibility of $L$-factors with
Langlands classification (\cite{Sha90b}). For $p$-adic fields
this is Corollary 8.2 of \cite{Sha92} in the
square-integrable case and follows from multiplicativity
(\cite{Sha90a}) in the general case. Note that on the
left-hand side we may take the epsilon factor as defined by Jacquet,
Piatetski-Shapiro and Shalika (\cite{JP-SS83},
\cite{JS90b}); it coincides with the one defined by Shahidi;
 see \cite{Sha84}. To finish the proof of Theorem
\ref{Main2} it remains to note that
$\eps(\frac12,\pi\otimes\dual{\pi})=1$ for any cuspidal
representation $\pi$ of ${\rm GL}_n(\A)$. This follows at once from the
next lemma which, at least in the nonarchimedean case, was
proved (even without the genericity assumption) by Bushnell and
Henniart (\cite{BH99}).
\enddemo

\specialnumber{4}\proclaim{Lemma} \label{rankin}
For any generic representation $\pi_v$ of ${\rm GL}_n$ over a local
field of characteristic $0${\rm ,}
\begin{equation} \label{eps}
\eps\left(\frac12,\pi_v\otimes\dual{\pi_v},\psi_v\right)=\omega_{\pi_v}(-1)^
{n-1},
\end{equation}
where $\omega_{\pi_v}$ is the central character of $\pi_v$.
\endproclaim

\section{Local analysis} \label{locanal}
\advance\eqcount by 7

In this section we prove Lemmas \ref{local}--\ref{rankin} which
were left out in the discussion of the previous section.

For the rest of the paper let $F$ be a local field of
characteristic $0$. We will suppress the subscript $v$ from all
notation and fix a nontrivial character $\psi$ of $F$ throughout.
As before, the $F$-points of an algebraic group $X$ over $F$ will
often  be denoted by $X$. We denote by $\nu$ the absolute value of
the determinant, viewed as a character on any one of the groups
${\rm GL}_n(F)$. If $\pi$ is a representation of ${\rm GL}_n$ and
$s\in\C$ we let $\pi\nu^s$ be the representation obtained   by
twisting $\pi$ by the character $\nu^s$. Let $\Irr_n$ be the set
of equivalence classes of irreducible (admissible) representations
of ${\rm GL}_n$. Given representations $\pi_i$, $i=1,\dots,k$ of
${\rm GL}_{n_i}$ we denote by $\pi_1\times\dots\times\pi_k$ the
representation on ${\rm GL}_n$ with $n=n_1+\dots+n_k$ induced from
the representation $\pi_1\otimes\dots\otimes\pi_k$ on the
parabolic subgroup of ${\rm GL}_n$ of type $(n_1,\dots,n_k)$.

\demo{{\rm 3.1.} Proof of Lemma \ref{rankin}} For completeness we
include a proof which was communicated to us by Herv\'e Jacquet.
We are very grateful to him.

By the functional equation the left-hand side of (\ref{eps}) is
$\pm1$. We prove the lemma by induction on $n$. If $\pi$ is not
essentially square-integrable then we can write
$\pi=\pi_1\times\pi_2$ where $\pi_i\in\Irr_{n_i}$ are generic.
We have
\begin{eqnarray*}
\hskip-36pt
&&\eps\left(\frac12,\pi_1\otimes\dual{\pi_2},\psi\right)
\eps\left(\frac12,\pi_2\otimes\dual{\pi_1},\psi\right)\\
\hskip-36pt&&\hskip1in =
\eps\left(\frac12,\pi_1\otimes\dual{\pi_2},\psi\right)
\eps\left(\frac12,\pi_2\otimes\dual{\pi_1},\conj{\psi}\right)\omega_{\pi_1
}^{n_2}(-1)
\omega_{\pi_2}^{n_1}(-1)\\
\hskip-36pt&&\hskip1in = \omega_{\pi_1}^{n_2}(-1)
\omega_{\pi_2}^{n_1}(-1)
\end{eqnarray*}
by the functional equation (\cite[p. 396]{JP-SS83}) and the
dependence of epsilon on $\psi$. By ``multiplicativity'' of
epsilon factors (loc.\ cit., p. 452) we get
$$
\eps\left(\frac12,\pi\otimes\dual{\pi},\psi\right)=
\eps\left(\frac12,\pi_1\otimes\dual{\pi_1},\psi\right)
\eps\left(\frac12,\pi_2\otimes\dual{\pi_2},\psi\right)
\omega_{\pi_1}^{n_2}(-1)\omega_{\pi_2}^{n_1}(-1)
$$
and we may use the induction hypothesis. Thus, it remains to
consider the case where $\pi$ is essentially square-integrable,
which immediately reduces to the case where $\pi$ is
square-integrable. In this case the zeta integral
$$
\Psi(s,W,W',\Phi)=\int_{N_n\bs {\rm
GL}_n}W(g)W'(g)\Phi((0,\dots,0,1)g)\abs{\det{g}}^s\ dg
$$
converges for $\Re(s)>0$ (loc.\ cit., (8.3)). Here $W$, $W'$ are
elements in the Whittaker spaces of $\pi$ and $\dual{\pi}$
respectively, and $\Phi$ is a Schwartz function on $F^n$. In
particular, $L(s,\pi\otimes\dual{\pi})$ has no pole (or zero) for
$\Re(s)>0$ and by the local functional equation (loc.\ cit., p.
391) we get
\begin{equation} \label{eps2}
\Psi\left(\frac12,\dual{W},\conj{\dual{W}},\hat\Phi\right)=\eps\left(\frac12
,\pi\otimes\dual{\pi},\psi\right)
\omega_{\pi}(-1)^{n-1}\Psi\left(\frac12,W,\conj{W},\Phi\right)
\end{equation}
for any $W$ and $\Phi$. Choose $W\not\equiv0$ and let $g$ be such
that $W(g)\ne0$. We may choose $\Phi\ge0$ such that
$\Phi((0,\dots,0,1)g)\ne0$ and $\hat{\Phi}\ge0$. For example, we
may take $\Phi$ of the form $\Phi_1\conv\Phi_1^\vee$ where
$\Phi_1\ge0$. Then clearly, both zeta integrals in (\ref{eps2})
are nonnegative and the one on the right-hand side is nonzero.
Hence $\eps(\frac12,\pi\otimes\dual{\pi},\psi)$ has the same sign
as $\omega_{\pi}(-1)^{n-1}$ and consequently, it is equal to it.
This finishes the proof of Lemma \ref{rankin}
\enddemo

If $\pi\in\Irr_n$ we denote by $\expo(\pi)$ the (central)
\emph{exponent} of $\pi$. It is the unique real number so that
$\pi\nu^{-\expo(\pi)}$ has a unitary central character. If
$\pi_1$, $\pi_2$ are generic and irreducible we let
$\intgln(\pi_1,\pi_2)$ be the normalized intertwining operator
$\pi_1\times\pi_2\rightarrow\pi_2\times\pi_1$ (depending on
$\psi$) as defined by Shahidi (\cite{Sha90b}) provided that it is
holomorphic there.

We recall that if $\pi$ and $\pi'$ are essentially square\/{\rm
-}\/integrable and $\abs{\expo(\pi)-\expo(\pi')}\break<1$ then
$\pi\times\pi'$ is irreducible and
$\pi\times\pi'\isom\pi'\times\pi$.

Recall the classification of the irreducible generic unitarizable
representations of ${\rm GL}_n$. (This is a very special case of
\cite{Tad86} in the $p$-adic case and \cite{Vog86} in the
archimedean case; cf.\ \cite{JS81} for the unramified case.) These
are the representations of the form
\begin{equation} \label{igu}
\sigma_1\times\dots\times\sigma_s\times\tau_1\nu^{\gamma_1}\times
\tau_1\nu^{-\gamma_1}\times\dots\times\tau_t\nu^{\gamma_t}\times\tau_t\nu^{-
\gamma_t}
\end{equation}
where the $\sigma_i$'s and the $\tau_j$'s are square integrable
(unitary), the $\sigma_i$'s are mutually inequivalent and $0\le
\gamma_j<\frac12$. Moreover, the data
$(\sigma_i)_{i=1}^s,(\tau_j,\gamma_j)_{j=1}^t$ are uniquely
determined up to permutation. Clearly, $\pi$ is self-dual if and
only if
$\{\sigma_i,\tau_j\nu^{\gamma_j}\}=\{\dual{\sigma_i},\dual{\tau_j}\nu^{\gamma_j}\}$
as multi-sets. Let $\Pi^{\sdu}$ be the set of self-dual generic
irreducible unitarizable representations of ${\rm GL}_n$.

Let $\SSS=\{S_n\}_{n\ge0}$ be any one of the families $\BBB={\rm
SO}(2n+1)$, $\CCC={\rm Sp}_n$ or $\DDD={\rm SO}(2n)$ (with
$S_0=1$). The family will be fixed throughout. In each case,
except for ${\rm SO}(2)$, the group $G=S_n$ is semisimple of rank
$n$ and we enumerate its simple roots
$\{\alpha_1,\dots,\alpha_n\}$ in the standard way. Recall the
automorphisms $\theta$ and $\Outer$ of $G$ defined in the previous
section. If $\pi$ is a representation of $G$ we let $\theta(\pi)$
be the representation obtained by twisting by $\theta$. Similarly
for $\Outer(\pi)$. We let $\Irr(S_n)$ be the set of equivalence
classes of irreducible representations of $S_n$. Let $\pi_i$,
$i=1,\dots,k$, be representations of ${\rm GL}_{n_i}$ and $\sigma$
a representation of $S_m$. Let $n=n_1+\dots+n_k+m$ and $Q$ be the
parabolic subgroup of $S_n$ obtained by ``deleting'' the simple
roots
$\alpha_{n_1},\alpha_{n_1+n_2},\dots,\alpha_{n_1+\dots+n_k}$, as
well as $\alpha_n$ in the case where $\SSS=\DDD$ and $m=1$. The
Levi subgroup $L$ of $Q$ is isomorphic to ${\rm
GL}_{n_1}\times\dots\times {\rm GL}_{n_k}\times S_m$.
As in \cite{Tad98} we denote by
$\pi_1\times\dots\times\pi_k\sdp\sigma$ the representation of
$S_n$ induced from the representation
$\pi_1\otimes\dots\otimes\pi_k\otimes\sigma$ of $Q$. We have,
$\pi\times\tau\sdp\sigma=\pi\sdp(\tau\sdp\sigma)$. In the case
$\SSS=\DDD$ we have $\Outer(\pi\sdp\sigma)=\pi\sdp\Outer(\sigma)$
for $\pi\in\Irr_n$ and $\sigma\in\Irr(S_m)$ with $m\ge1$.

Let $L$ be a Levi subgroup of $G$ and let $w_0$ (resp. $w_0^L$) be
the longest element in the Weyl group of $G$ (resp. $L$). We
denote by $w_L$ the Weyl group element $w_0w_0^L$. In particular
$w_M$ is defined, where we recall that $M\isom {\rm GL}_n$ is the Siegel
Levi.

Suppose that $\pi_i\in\Irr_{n_i}$ are essentially square
integrable with $\expo(\pi_1)>\expo(\pi_2)>\dots>\expo(\pi_k)>0$
and that $\sigma\in\Irr(S_m)$ is square integrable. Let $Q$ and
$L$ be as before. Then
\begin{itemize}
\item[1.]   $\pi_1\times\dots\times\pi_k\sdp\sigma$ admits a unique
irreducible quotient.
\item[2.] The multiplicity of this quotient in the semi-simplification of
$\pi_1\times\dots\times\pi_k\sdp\sigma$ is one.
\item[3.] The quotient is isomorphic to the image of the (unnormalized)
intertwining
operator
$$
M_w:\pi_1\times\dots\times\pi_k\sdp\sigma\rightarrow\Outer^{n_1+\dots+n_k}
(\opp{\pi_1}\times\dots\times\opp{\pi_k}\sdp\sigma)
$$
with respect to $Q$ and $w$ where $w=w_L$.
\item[4.] $M_w$ is given by a convergent integral.
\end{itemize}
This is the Langlands quotient in this setup. For all this see
\cite{BW00}. Let $Q'$ be the parabolic subgroup with Levi
subgroup $L'$ isomorphic to ${\rm GL}_{n_1+\dots+n_k}\times S_m$ and let
$\pi=\pi_1\times\dots\times\pi_k$. The operator $M_w$ is obtained
as the composition of the intertwining operator
\begin{equation} \label{inter}
\pi\sdp\sigma\rightarrow
\Outer^{n_1+\dots+n_k}(\opp{\pi}\sdp\sigma)
\end{equation}
with respect to $Q'$ and $w_{L'}$, and an intertwining operator
$M_2$ ``inside''\break ${\rm GL}_{n_1+\dots+n_k}$. Under the weaker
hypothesis
that $\expo(\pi_1)\ge\dots\ge\expo(\pi_k)>0$ the statements
1--3 will continue to hold provided that
$M_2$ is normalized. This is because the $R$-groups for general
linear groups are trivial. In particular, if $\pi$ is irreducible
then the Langlands quotient is isomorphic to the image of the
intertwining operator (\ref{inter}).

If $\pi$ is a representation of ${\rm GL}_n$ we let
$I(\pi,s)=I^G(\pi,s)$ be the induced representation
$\pi\nu^s\sdp\triv$. Similar notation will be used for induction
from the parabolic subgroup $\theta(P)$. We denote by
$\Mabs(\pi,s)=\Mabs(s):I(\pi,s)\rightarrow
I(\w{\pi},-s)=\theta(I(\opp{\pi},-s))$ the unnormalized
intertwining operator with respect to $P$ and $w_M$. If $\pi$ is
generic we denote by $\nrmint(\pi,s)=\nrmint(s)$ the normalized
intertwining operator (with respect to $P$ and $w_M$). In the case
$G={\rm SO}(2)$ we set $\Mabs(s)=\nrmint(s)=\id$.

We will often use the following fact. Suppose that
$\pi=\pi_1\times\dots\times\pi_k$ is a generic representation of
${\rm GL}_n$ with $\pi_i\in\Irr_{n_i}$. We may identify $I(\pi,s)$ with
$\Ind_Q^G(\pi_1\nu^s\otimes\dots\otimes\pi_k\nu^s)$ where
$\pi_1\otimes\dots\otimes\pi_k$ is viewed as a representation of
the parabolic subgroup $Q$ of $G$ whose Levi subgroup is the Levi
subgroup of ${\rm GL}_n$ of type $(n_1,\dots,n_k)$. We may also identify
$\opp{\pi}$ with $\opp{\pi_k}\times\dots\times\opp{\pi_1}$ and
$I(\opp{\pi},s)$ with
$\Ind_{\opp{Q}}^G\opp{\pi_k}\nu^s\otimes\dots\otimes\opp{\pi_1}\nu^s$.
Under these identifications $\nrmint(s)$ becomes the normalized
intertwining operator
$$
\Ind_Q^G(\pi_1\nu^s\otimes\dots\otimes\pi_k\nu^s)\rightarrow
\theta(\Ind_\opp{Q}^G(\opp{\pi_k}\nu^{-s}\otimes\dots\otimes\opp{\pi_1}\nu^{
-s}))
$$
with respect to $Q$ and $w_M$. This is merely a reformulation of
the multiplicativity of $L$ and $\eps$-factors
(\cite{Sha90a}). As a result, we may decompose the operator
$\nrmint(\pi,s)$ as a product of ``basic'' intertwining operators
according to the reduced decomposition of $w_M$. Each basic
intertwining operator is obtained by inducing an operator of the
form $\nrmint(\pi_i,s)$ or $\intgln(\pi_i,\pi_j)$ or
$\intgln(\pi_i,\opp{\pi_j})$ with $i>j$.
\enddemo

 3.2. {\it Proof of Lemma} \ref{local}. Let $\pi\in\Pi^\sdu$
and $\Re(s)\ge\frac12$. We may write $\pi\nu^s$ as
$\pi_1\times\dots\times\pi_k$ with $\pi_i$ essentially
square-integrable and $\expo(\pi_1)\ge\dots\ge\expo(\pi_k)>0$.
Hence $I(\pi,s)$ admits a Langlands quotient, which by the
discussion above, is given by the image of $\Mabs(\pi,s)$. By
multiplicativity, the statements about the $L$-functions follow
from the holomorphy of $L(s,\pi_i)$, $L(s,\pi_i,\sym^2)$,
$L(s,\pi_i,\wedge^2)$ and $L(s,\pi_i\otimes\pi_j)$ at $s=0$, which
in turn follows from
 \cite[Prop.~7.2]{Sha90b}. The statements about the
normalizing factors and $\nrmint(\pi,s)$ follow immediately.
\hfill\qed\vglue12pt

\demo{{\rm 3.3.} Proof of Lemma \ref{claim}} Recall the definition
of the operators $\herm^\psi(\pi,s)$. (We assume that $\theta=\id$
and that $\pi\in\Irr_n$ is self-dual.) We first note that
$\ident_\pi$ is Hermitian since, being an intertwining operator of
order two, it must preserve the inner product. We conclude that
$\ident_{\pi,s}^*=\ident_{\w\pi,-\overline{s}}$ where $^*$ denotes
the Hermitian dual. Also, a direct calculation shows the relation
$$
\Mabs(\w\pi,s)\ident_{\w\pi,s}=\iota_{\pi,-s}\Mabs(\pi,s).
$$
Moreover, by (\ref{w2}) the Hermitian dual of $\Mabs(\pi,s)$ is given by
$\omega_{\pi}(-1)^k\Mabs(\w\pi,\overline{s})$ where $k=n$ if
$G$ is symplectic and $k=n+1$ if $G$ is orthogonal.
On the other hand, by the dependence of root numbers on
the additive character it is easily deduced that
$$
\overline{\m^\psi(\pi,s)}=\omega_\pi(-1)^k\m^\psi(\pi,\overline{s}).
$$
The Hermitian property of $\herm^\psi(\pi,s)$ for $s$ real follows.

To prove the second part we use the argument
of \cite[Prop.~6.3]{KS88}. Let $W_\pi^\psi(\cdot,s)$ be the
Whittaker functional on $I(\pi,s)$ and let
$W_{\w\pi}^\psi(\cdot,s)$ be the Whittaker functional on
$I(\w\pi,s)$ obtained through $\iota_\pi$. They are holomorphic,
nonzero (\cite{Sha81}), and satisfy the functional equation
\begin{equation} \label{funcequ}
W_\pi^\psi(\varphi,s)=
c(\pi,s,\psi)W_{\w\pi}^\psi(M(\pi,s)\varphi,-s)
\end{equation}
where $c^\psi(\pi,s)$ is the ``local coefficient'' which was studied by
Shahidi. By \cite[Th.~3.5]{Sha91} it is given by
\begin{equation} \label{plancherel}
c^\psi(\pi,s)= \left\{ \begin{array}{ll}
\frac{\eps(2s,\pi,r,\psi^{-1})L(1-2s,\pi,r)}{L(2s,\pi,r)}&\SSS=\BBB,
\DDD\\[5pt]
\frac{\eps(s,\pi,\psi^{-1})L(1-s,\pi)}{L(s,\pi)}\cdot
\frac{\eps(2s,\pi,r,\psi^{-1})L(1-2s,\pi,r)}{L(2s,\pi,r)}&\SSS=\CCC
\end{array}\right.
\end{equation}
where $r=\sym^2$ for $\SSS=\BBB$ and $r=\wedge^2$ for $\SSS=\CCC,
\DDD$. (Here we use that $\pi$ is self-dual.) By the identification
$\ident_\pi:\w\pi\rightarrow\pi$, (\ref{funcequ}) becomes
$$
W_\pi^\psi(\varphi,s)=c^\psi(\pi,s)\m^\psi(\pi,s)
W_\pi^\psi(\herm^\psi(\pi,s)\varphi,-s).
$$
The term $c^\psi(\pi,s)m^\psi(\pi,s)$ is either
$L(1-2s,\pi,r)/L(1+2s,\pi,r)$ if $\SSS=\BBB, \DDD$ or
$L(1-2s,\pi,r)/L(1+2s,\pi,r)\cdot L(1-s,\pi)/L(1+s,\pi)$ if
$\SSS=\CCC$. It follows that
$$
\herm^\psi(\pi,-s)\herm^\psi(\pi,s)=I.
$$
We infer that $\herm^\psi(\pi,s)$ is unitary, and in particular,
holomorphic at $s=0$. Moreover, $\herm^\psi(\pi,0)$ fixes the
$\psi$-generic irreducible constituent of $I(\pi,0)$, since
$L(s,\pi)$ and $L(s,\pi,r)$ are holomorphic at $s=1$ by Lemma
\ref{local}.
\enddemo

\vglue3pt

The rest of the paper is devoted to the proof of Lemma
\ref{claim2}. Since the lemma is evidently independent of the
choice of the character $\psi$, we will suppress it from the
notation.

\demo{{\rm 3.4.} Representations of $G$\/{\rm -}\/type} Let
$\sigma$ be a self-dual square-integrable representation of ${\rm
GL}_n$ and suppose that $\theta=\id$. By the theory of $R$-groups
(e.g.\ \cite{Gol94}) the following conditions are equivalent.
\begin{itemize}
\item[1.] $I(\sigma,0)$ is irreducible.
\item[2.] $\herm(\sigma,0)$ is a scalar.
\item[3.] The Plancherel measure $\mu(\sigma,s)$ is zero at $s=0$.
\end{itemize}

\demo{Definition {\rm 1}}
An essentially square-integrable representation $\sigma$ of ${\rm GL}_n$
will be called of $G$-{\it type} (or of $\SSS$-{\it type} if we do not want
to
specify $n$) if it is self-dual (in particular,
$\expo(\sigma)=0$), $\theta=\id$, and the conditions above are
satisfied.
\enddemo

\specialnumber{1}\proclaim{Proposition} \label{char}
Let $\sigma$ be a square\/{\rm -}\/integrable representation of ${\rm
GL}_n$. Then
$I^G(\sigma,s)$ is irreducible for $0<s<1$ except possibly for
$s=\frac12$. Moreover{\rm ,} if $I(\sigma,\frac12)$ is reducible then
$\sigma$ is of $G$\/{\rm -}\/type.
\endproclaim

\demo{Proof} By the results of Muic (\cite{Mui01}) we may use
Proposition 5.3 of \cite{CS98}. Thus the reducibility points of
$I(\sigma,s)$ for $s>0$ are the poles of $L(1-2s,\sigma,r)$ (if
$\SSS=\BBB$ or $\DDD$) or $L(1-s,\sigma)L(1-2s,\sigma,r)$ (if
$\SSS=\CCC$). These $L$-functions are computed in
\cite[Prop.~8.1]{Sha92}. In particular, $L(s,\sigma,r)$ is
holomorphic for $s>-1$ except possibly for $s=0$ and $L(s,\sigma)$
is holomorphic for $s>0$. Therefore $I(\sigma,s)$ is irreducible
for $0<s<1$ except possibly for $s=\frac12$ and moreover, if
$I(\sigma,\frac12)$ is reducible then $L(s,\sigma,r)$ has a pole
at $s=0$. In the latter case $\theta=\id$, $\sigma$ is self-dual
and the local coefficient vanishes at $0$ (loc.\ cit.). By
\cite[(1.4)]{Sha90b} the same will be true for the Plancherel
measure. \phantom{anytime}
\enddemo

{\it Remark} 1.
Shahidi also proved the following in (\cite{Sha92}). Suppose
that $\sigma$ is a self-dual square-integrable representation of
${\rm GL}_n$ which is not the trivial character of ${\rm GL}_1$. Then the
following are equivalent:
\begin{itemize}
\item[1.] $\sigma$ is of ${\rm Sp}_n$ type.
\item[2.] $\sigma$ is not of ${\rm SO}(2n+1)$ type.
\item[3.] $\sigma$ is of ${\rm SO}(2n)$ type.
\end{itemize}
In particular, in this case $n$ must be even. We will not use this
fact.

For convenience, we consider the set $\Pi^{\sd}$ of all
representations of the form
\begin{equation} \label{sd}
\sigma_1\times\dots\times\sigma_s\times\tau_1\times
\opp{\tau_1}\times\dots\times\tau_t\times\opp{\tau_t}
\end{equation}
where the $\sigma_i$'s are square-integrable, self-dual and (as we
may assume) mutually inequivalent, and the $\tau_j$'s are
essentially square-integrable with $0\le\expo(\tau_j)<\frac12$.
Any element of $\Pi^\sd$ is irreducible, generic and self-dual.
Clearly, $\Pi^\sd\supset\Pi^\sdu$. The condition on
$\pi\in\Pi^\sd$ to belong to $\Pi^\sdu$ (i.e. to be unitarizable)
is that each $\tau_j$ which is not square-integrable appears in
(\ref{sd}) the same number of times as
$\opp{\tau_j}\nu^{2\expo(\tau_j)}$. If $\pi\in\Pi^\sd$ then by the
discussion of subsection~3.2, $I(\pi,\frac12)$ admits a Langlands
quotient, which will be denoted by $\LQ{\pi}$. It is obtained as
the image of $\nrmint(\pi,\frac12)$ (or $\Mabs(\pi,\frac12)$).

Also, if $\chi$ is an essentially square-integrable representation
of ${\rm GL}_n$ we denote by $\lqq{\chi}$ the unique irreducible
quotient of $\chi\nu^\frac12\times\chi\nu^{-\frac12}$. It is
isomorphic to the image of the intertwining operator
$\intgln(\chi\nu^\frac12,\chi\nu^{-\frac12})$.

\specialnumber{5}\proclaim{Lemma} \label{simple} Let $\chi$ be an
essentially square\/{\rm -}\/integrable representation of ${\rm
GL}_n$ with $0\le\expo(\chi)<\frac12$. Assume that $\chi$ is not
of $\SSS$\/{\rm -}\/type. Then
$\LQ{\chi\times\opp{\chi}}\isom\Outer^n(\lqq{\chi}\sdp\triv)$.
\endproclaim

\demo{Proof}
The Langlands quotient is obtained as the image of the longest
intertwining operator, which is the composition of the following
intertwining operators:
$$
I\left(\chi\times\opp{\chi},\frac12\right)\stackrel{\id\sdp\nrmint(\opp{\chi},\frac12)}{\lrar}
\chi\nu^\frac12\sdp\Outer^n\left(\chi\nu^{-\frac12}\sdp\triv\right)\isom
\Outer^n\left(\chi\nu^\frac12\times\chi\nu^{-\frac12}\sdp\triv\right)$$
$$\stackrel{\Outer^n(
I(R_1,0))}{\lrar}
\Outer^n\left(\chi\nu^{-\frac12}\times\chi\nu^\frac12\sdp\triv\right)\isom
\chi\nu^{-\frac12}\sdp\Outer^n\left(\chi\nu^\frac12\sdp\triv\right)$$
$$
\stackrel{\id\sdp\Outer^n (\nrmint (\chi,\frac12 ) )}{\lrar}
\chi\nu^{-\frac12}\times\opp{\chi}\nu^{-\frac12}\sdp\triv
$$
where $R_1=\intgln(\chi\nu^\frac12,\chi\nu^{-\frac12})$. By
Proposition \ref{char} the only map which is not an isomorphism is
$\Outer^n( I(R_1,0))$, whose image is
$\Outer^n(\lqq{\chi}\sdp\triv)$ as required.
\enddemo

Any $\pi\in\Pi^\sd$ can be written uniquely as
$\pi^\orth\times\pi^\orthpairs\times\pi^\symp$ with
\begin{itemize}
\item $\pi^\symp$ of the form $\sigma_1\times\dots\times\sigma_s$
where the $\sigma_i$'s are square-integrable, self-dual and of
$\SSS$-type;
\item $\pi^\orth$ of the form $\rho_1\times\dots\times\rho_r$
where the $\rho_i$'s are square-integrable, mutually inequivalent,
self-dual and not of $\SSS$-type;
\item $\pi^\orthpairs$ of the form
$\tau_1\times\opp{\tau_1}\times\dots\times\tau_t\times\opp{\tau_t}$
where the $\tau_j$'s are essentially square-integrable, not of
$\SSS$-type (self-dual or not), and $0\le\expo(\tau_j)<\frac12$.
\end{itemize}
Note that $\pi^\orth$ and $\pi^\symp$ are tempered.

\demo{Definition {\rm 2}}
We say that $\pi\in\Pi^{\sd}$ is of $G$-type if $\pi^\orth=0$.
\enddemo

The definition is suggested by the local Langlands reciprocity.
Note that if $\pi$ is of ${\rm SO}(2n)$ type then $n$ is even.

The crucial property of representations of $G$-type is the following.

\specialnumber{6}\proclaim{Lemma} \label{imply} If $\pi\in\Pi^\sd$
is of $G$\/{\rm -}\/type then $\herm(0)$ is a nonzero scalar.
\endproclaim

\demo{Proof}
We use induction on $n$, the case $n=0$ being trivial. For the
induction step, we can assume that $\pi=\pi'\times\omega$ where
$\pi'\in\Pi^\sd$ is of $\SSS$-type and $\omega\in\Irr_l$ is either
square-integrable and of $\SSS$-type or of the form
$\tau\times\opp{\tau}$ where $\tau\in\Irr_m$ is essentially
square-integrable. Note that $l$ is even if $\SSS=\DDD$. The
operator $\nrmint(0)$ can be written as the composition of the
following intertwining operators:
\begin{eqnarray} \label{4step}
&&I(\pi'\times\omega,0)\stackrel{\id\sdp\nrmint(\omega,0)}{\lrar}
I(\pi'\times\opp{\omega},0)\stackrel{I(\nrmint_1,0)}{\lrar}\\
&&\qquad\qquad
I(\opp{\omega}\times\pi',0)\stackrel{\id\sdp\nrmint(\pi',0)}{\lrar}
I(\opp{\omega}\times\opp{\pi'},0)\isom
I(\opp{(\pi'\times\omega)},0)\nonumber
\end{eqnarray}
where $\nrmint_1=\intgln(\pi',\opp{\omega})$. The last
identification is induced by the isomorphism
$\opp{\omega}\times\opp{\pi'}\isom\opp{(\pi'\times\omega)}$. By
the induction hypothesis the third map is a nonzero scalar
multiple of $\id\sdp\ident(\pi',0)^{-1}$. Also, by uniqueness,
$\ident_\pi:\opp{\omega}\times\opp{\pi'}\rightarrow\pi$ is a
scalar multiple of
$(\id\times\ident_\omega)\nrmint_1^{-1}(\id\times\ident_\pi')$.
All in all, the map $\herm(\pi,0)$ is a scalar multiple of
$$
I((\id\times\ident_\omega)\nrmint_1^{-1}(\id\times\ident_\pi'),0)
\circ\id\sdp\ident(\pi',0)^{-1}\circ I(\nrmint_1,0)\circ
\id\sdp\nrmint(\omega,0) =\id\sdp\herm(\omega,0).
$$
It remains to show that $\herm(\omega,0)$ is a scalar in the two
cases above. In the first case, this follows from the definition
of $\SSS$-type. In the second case, we decompose
$\nrmint(\omega,0)$ as before as
\begin{eqnarray*}
&&I(\tau\times\opp{\tau},0)\xrightarrow{\id\sdp\nrmint(\opp{\tau},0)}
\Outer^m(I(\tau\times\tau,0))\xrightarrow{\Outer^m(I(\nrmint_2,0))}\\
&&\qquad
\Outer^m(I(\tau\times\tau,0))\xrightarrow{\id\sdp\Outer^m(\nrmint(\tau,0))}
I(\tau\times\opp{\tau},0)\isom I(\opp{(\tau\times\opp{\tau})},0).
\end{eqnarray*}
Note that the map $\nrmint_2=\intgln(\tau,\tau)$ is a scalar, and
similarly for the map
$\ident_\omega:\tau\times\opp{\tau}\isom\opp{(\tau\times\opp{\tau})}\rightarrow
\tau\times\opp{\tau}$. Thus, $\herm(\omega,0)$ is a scalar
multiple of
$\id\sdp\Outer^m(\nrmint(\tau,0))\circ\nrmint(\opp{\tau},0)$ which
is $\id$ by the properties of the normalized intertwining
operator. \phantom{takeawalk}
\enddemo

\vglue-8pt
{\it Remark} 2.
The converse to Lemma \ref{imply} is also true.

\vglue6pt  3.5. {\it Langlands quotient}.

We extract a few results from \cite{MW89} (cf.\ I.6.3 for the
$p$-adic case and I.7 for the archimedean case).
\enddemo

\specialnumber{7}\proclaim{Lemma} \label{GLn} Let $\pi$ and $\pi'$
be irreducible representations of ${\rm GL}_n$ and ${\rm GL}_{n'}$
respectively.
\begin{itemize}
\ritem{1.} If $\pi\times\pi'$ is irreducible then
$\pi\times\pi'\isom\pi'\times\pi$. \ritem{2.}   Let $\pi$ and
$\pi'$ be essentially square\/{\rm -}\/integrable. Suppose that
$\abs{\expo(\pi)-\expo(\pi')}\break<1$. Then $\pi\times\pi'$,
$\pi\nu^\frac12\times\lqq{\pi'}$ and $\lqq{\pi}\times\lqq{\pi'}$
are irreducible. \ritem{3.}   Suppose that $\pi$ and $\pi'$ are
inequivalent square\/{\rm -}\/integrable representations. Then
$\pi\nu^\gamma\times\pi'\nu^{-\frac12}$ is irreducible for
$-1<\gamma<1$.
\end{itemize}

\endproclaim

We will also need the following lemma which is based on
\cite{Jan96}.

\specialnumber{8}\proclaim{Lemma} \label{spn} Let
$\pi_i\in\Irr_{n_i}$ for $i=1,\dots,k$ and $\sigma\in\Irr(S_m)$.
Suppose that $\pi_i\times\pi_j${\rm ,} $\pi_i\times\opp{\pi_j}$
are irreducible for all $i\ne j$ and $\pi_i\sdp\sigma$ is
irreducible for all $i$. Then
\begin{equation}
\pi_1\times\dots\times\pi_k\sdp\sigma\isom\Outer^{n_1+\dots+n_k}
(\opp{\pi_1}\times\dots\times\opp{\pi_k}\sdp\sigma).
\end{equation}
Suppose in addition that the $\pi_i$\/{\rm '}\/s are essentially
square\/{\rm -}\/integrable with $\expo(\pi_i)>0$ and $\sigma$ is
square\/{\rm -}\/integrable. Then
$\pi_1\times\dots\times\pi_k\sdp\sigma$ is irreducible.
\endproclaim

\demo{Proof} In the case where $k=1$ we note that if
$\pi\in\Irr_n$ and $\sigma\in\Irr(S_m)$ then
$\pi\sdp\sigma=\Outer^n(\opp{\pi}\sdp\sigma)$ in the Grothendieck
group since $\pi\otimes\sigma$ and
$\Outer^n(\opp{\pi}\otimes\sigma)$ are associate.
The case $k>1$ and the last statement
are proved in (\cite{Jan96}) for the cases $\SSS=\BBB, \CCC$. The
proof carries over almost literally (except for putting in some
$\Outer$'s) to the case $\SSS=\DDD$ (cf.\ Proposition \ref{lqsetup}
below).
\enddemo

Let $\pi\in\Pi^\sd$. Recall that $I(\pi,\frac12)$ admits a
Langlands quotient, denoted by $\LQ{\pi}$, which is isomorphic to
the image under $\nrmint(\pi,\frac12)$.
\specialnumber{2}\proclaim{Proposition} \label{lqsetup} Let
$\pi=\pi^\orth\times\pi^\orthpairs\times\pi^\symp\in\Pi^\sd$ be as
above. Then
\begin{equation} \label{lq}
\LQ{\pi}\isom
\Outer^\eps(\lqq{\tau_1}\times\dots\times\lqq{\tau_t}\times\pi^\orth\nu^\frac12\sdp
\LQ{\pi^{\symp}})\qquad
\end{equation}
for $\eps$ either $0$ or $1$ (depending only on $\pi^\orthpairs$).
Hence{\rm ,}
$$
\LQ{\pi}\isom\pi^\orth\nu^\frac12\sdp\LQ{\pi^\orthpairs\times\pi^\symp}.
$$
\endproclaim

\vglue-12pt
{\it Proof}.
Clearly, the second statement follows from the first. Let
$\Lambda$ be the right-hand side of (\ref{lq}). Following the
argument of \cite[Th.~3.3]{Jan96} we will argue that
\begin{eqnarray}
&&\Lambda\hbox{ is a quotient of }I(\pi,{\textstyle\frac12})\hbox{ for
$\eps$ either $0$ or $1$}.\label{quotient}\\
&&\hskip1in \Lambda\hbox{ is irreducible}.\label{irred}
\end{eqnarray}
The first statement is proved by induction on $n$, as in
subsection~3.4. Since the case where $\pi^\orthpairs=0$ is
immediate, we may assume for the induction step that
$\pi=\pi'\times\tau\times\opp{\tau}$ where $\pi'\in\Pi^\sd$,
$\tau\in\Irr_m$ is essentially square-integrable,
$0\le\expo(\tau)<\frac12$ and $\tau$ is not of $\SSS$-type.

It follows from Lemma \ref{simple} that up to $\Outer$,
$I(\pi,\frac12)$ has a quotient isomorphic to
$I(\pi'\nu^\frac12\times\lqq{\tau},0)$. It follows from part
2 of Lemma \ref{GLn} that
$\pi'\nu^\frac12\times\lqq{\tau}$ is irreducible, and hence, that
$\pi'\nu^\frac12\times\lqq{\tau}\isom\lqq{\tau}\times\pi'\nu^\frac12$.
We deduce that up to $\Outer$, $I(\pi,\frac12)$ admits
$I(\lqq{\tau}\times\pi'\nu^\frac12,0)$, and thus also
$\lqq{\tau}\sdp\LQ{\pi'}$, as a quotient. This implies
(\ref{quotient}) by the induction hypothesis.

To prove (\ref{irred}), it suffices to show that
$\theta(\Lambda)\isom\dual{\Lambda}$. (This condition does not
depend on $\eps$.) Indeed, we have
$\theta(\LQ{\pi})\isom\dual{\LQ{\pi}}^{\phantom{|}}\hskip-2pt$
(cf.\ \cite{Jan96}) since both sides are the unique irreducible
subrepresentation of $I(\dual{\pi},-\frac12)$ by (\ref{contra}).
We would conclude that $\LQ{\pi}$ is both a quotient and a
subrepresentation of $\Lambda$. However, $\LQ{\pi}$ is the unique
irreducible quotient of $\Lambda$, and it has multiplicity-one in
the semi-simplification of $\Lambda$. Thus,
$\Lambda\isom\LQ{\pi}$.

We shall write $\pi_3$ for $\pi^\symp$. To show that
$\theta(\Lambda)\isom\dual{\Lambda}$ we note once more that
$\dual{\LQ{\pi_3}}=\LQ{\pi_3}$. By Lemmas \ref{GLn} and \ref{spn}
it suffices to show that both $\rho\nu^\frac12\sdp\LQ{\pi_3}$ and
$\lqq{\tau}\sdp\LQ{\pi_3}$ are irreducible where $\rho\in\Irr_l$
is square-integrable self-dual not of $\SSS$-type and $\tau$ is as
before.

To prove this, consider the representation
$\pi'=\tau\times\opp{\tau}\times\pi_3$. The Langlands quotient of
$I(\pi',\frac12)$ is the image of the operator $M_{w_0}$ which is
the composition of the intertwining operators
\begin{eqnarray*}
&&\hskip-16pt
I\left(\tau\nu^\frac12\times\opp{\tau}\nu^\frac12\times\pi_3\nu^\frac12,0\right)
\\\longrightarrow\
&&\hskip-16pt
I\left(\tau\nu^\frac12\times\pi_3\nu^\frac12\times\opp{\tau}\nu^\frac12,0\right)
\\\longrightarrow\
&&\hskip-16pt
\Outer^m\left(I\left(\tau\nu^\frac12\times\pi_3\nu^\frac12\times\tau\nu^{-\frac12},0\right)\right)
\\\longrightarrow\
&&\hskip-16pt
\Outer^m\left(I\left(\tau\nu^\frac12\times\tau\nu^{-\frac12}\times\pi_3\nu^\frac12,0\right)\right)
\\\longrightarrow\
&&\hskip-16pt
\Outer^m\left(I\left(\tau\nu^{-\frac12}\times\tau\nu^\frac12\times\opp{\pi_3
}\nu^{-\frac12},0\right)\right)
\\\longrightarrow\
&&\hskip-16pt
\Outer^m\left(I\left(\tau\nu^{-\frac12}\times\opp{\pi_3}\nu^{-\frac12}\times
\tau\nu^\frac12,0\right)\right)
\\\longrightarrow\
&&\hskip-16pt
I\left(\tau\nu^{-\frac12}\times\opp{\pi_3}\nu^{-\frac12}\times\opp{\tau}\nu^
{-\frac12},0\right)
\\\longrightarrow\
&&\hskip-16pt
I\left(\tau\nu^{-\frac12}\times\opp{\tau}\nu^{-\frac12}\times\opp{\pi_3}\nu^
{-\frac12},0\right)
\\\longrightarrow\
&&\hskip-16pt
I\left(\opp{\tau}\nu^{-\frac12}\times\tau\nu^{-\frac12}\times\opp{\pi_3}\nu^
{-\frac12},0\right).
\end{eqnarray*}
Again by Lemma \ref{GLn} and Proposition \ref{char}, all arrows
except the fourth one are isomorphisms. Thus, the Langlands
quotient is isomorphic to the image of the fourth map, which is
$\Outer^m(\lqq{\tau}\sdp\LQ{\pi_3})$. Hence, the latter is
irreducible. Similarly, if $\pi'=\rho\times\pi_3$ then $\LQ{\pi'}$
is the image of the composition of the intertwining operators
\begin{eqnarray*}
&&\hskip-16pt I\left(\rho\nu^\frac12\times\pi_3\nu^\frac12,0\right)
\\\longrightarrow\
&&\hskip-16pt I\left(\rho\nu^\frac12\times\opp{\pi_3}\nu^{-\frac12},0\right)
\\\longrightarrow\
&&\hskip-16pt I\left(\opp{\pi_3}\nu^{-\frac12}\times\rho\nu^\frac12,0\right)
\\\longrightarrow\
&&\hskip-16pt
\Outer^l\left(I\left(\opp{\pi_3}\nu^{-\frac12}\times\opp{\rho}\nu^{-\frac12}
,0\right)\right).
\end{eqnarray*}
Again, all maps except the first are isomorphisms. Thus, as
before, $\rho\nu^\frac12\sdp\LQ{\pi_3}$ is irreducible.
\hfill\qed\vglue8pt

For future reference, let us reformulate the conclusion of
Proposition \ref{lqsetup}. Using a decomposition of $w_0$ we may
decompose $\nrmint(\pi,\frac12)$
as
\begin{eqnarray} && \label{seq}\\
\hskip-2pt&&\hskip-8pt I\left(\pi,\frac12\right)
=I\left(\xprod_{i=1}^r\left(\tau_i\nu^\frac12\times\opp{\tau_i}\nu^\frac12\right)\times
\pi^\orth\nu^\frac12\times\pi^\symp\nu^\frac12,0\right)\nonumber\\
\hskip-2pt&\srightarrow{R_{w_1}}&\hskip-8pt
\Outer^\eps\left(I\left(\xprod_{i=1}^r\left(\tau_i\nu^\frac12\times\tau_i\nu
^{-\frac12}\right)\times
\pi^\orth\nu^\frac12\times\pi^\symp\nu^\frac12,0\right)\right)\nonumber\\
\hskip-2pt&\srightarrow{R_{w_2}}&\hskip-8pt
\Outer^\eps\left(I\left(\xprod_{i=1}^r\left(\tau_i\nu^{-\frac12}\times\tau_i
\nu^\frac12\right)\times
\pi^\orth\nu^\frac12\times\opp{\pi^\symp}\nu^{-\frac12},0\right)\hskip-2pt\right)\nonumber\\
\hskip-2pt&\srightarrow{R_{w_3}}&\hskip-8pt
\Outer^{\eps'}\left(I\left(\xprod_{i=1}^r\left(\opp{\tau_i}\nu^{-\frac12}\times\tau_i\nu^{-\frac12}\right)
\times\opp{\pi^\orth}\nu^{-\frac12}\times\opp{\pi^\symp}\nu^{-\frac12},0\right)\hskip-2pt\right)\nonumber\\
\hskip-2pt&&\hskip-8pt  =
\theta\left(I\left(\pi,-\frac12\right)\hskip-2pt\right)\nonumber
\end{eqnarray}
where $R_{w_i}$ are normalized intertwining operators. We observe
that the image of $R_{w_2}$ (of the whole induced space) is
isomorphic to the right-hand side of (\ref{lq}), and hence it is
the Langlands quotient. By irreducibility and multiplicity-one of
Langlands quotient  $\im(R_{w_2}\circ R_{w_1})=\im(R_{w_2})$
and $\ker(R_{w_3}\circ R_{w_2})=\ker(R_{w_2})$.

\demo{{\rm 3.6.} Reduction to the tempered case} Let
$\pi\in\Pi^\sdu$. We may write $\pi=\pi^\temp\times\pi^\nontemp$
where $\pi^\temp\in\Pi^\sdu$ is tempered and $\pi^\nontemp$ is of
the form
$\xprod_i(\omega_i\nu^{\beta_i}\times\opp{\omega_i}\nu^{-\beta_i})$
with $\omega_i$ square-integrable and $0<\beta_i<\frac12$.
Clearly, $\pi^\nontemp$ appears as a factor of $\pi^\orthpairs$.
We will deform the nontempered parameters of $\pi$. For $0\le
t\le1$ let
$$
\pi_t=\pi^\temp\times\xprod_i\left(\omega_i\nu^{t\beta_i}\times
\opp{\omega_i}\nu^{-t\beta_i}\right)=\pi^\temp\times\pi_t^{\nontemp}.
$$
Then $\pi_t$ is a ``deformation'' in $\Pi^{\sdu}$ from
$\pi\isom\pi_1$ to the tempered representation $\pi_0$. Clearly
$\pi_t^\orth=\pi^\orth$ for all $t$ and $\pi_t^\symp=\pi^\symp$
for $t\ne0$ although not necessarily for $t=0$. The form
$\Inner(\pi,s)$ depends on the unitary structure on $\pi$, or what
amounts to the same, on a ${\rm GL}_n$-invariant positive-definite
Hermitian form on $\pi$. We identify the ambient vector spaces of
$\pi_t$ with that of $\pi$ in the usual way. The $K$-action does
not depend on $t$, where $K$ denotes the standard maximal compact.
We may choose a family of ${\rm GL}_n$-invariant positive-definite
Hermitian forms on $\pi_t$ which depends
continuously on $t$ (using intertwining operators for example). 

The following lemma will reduce Lemma \ref{claim2} to the tempered
case.
\enddemo

\specialnumber{9}\proclaim{Lemma} \label{reduce}
{\rm 1.} The definiteness of $\Inner(\pi_t,0)$ does not depend on $t$.
\vglue4pt {\rm 2.} If $\Inner(\pi,\frac12)$ is semi\/{\rm -}\/definite then
$\Inner(\pi_0,\frac12)$ is semi-definite with the same sign.
\endproclaim

We will use the following elementary lemma.

\specialnumber{10}\proclaim{Lemma} \label{elem} \hskip-8pt
Let $\{l_\beta\}_{a\le\beta\le b}$ be a continuous family of
Hermitian forms on~$\C^m$. Suppose that $\rank(l_\beta)$ is
constant for $a<\beta\le b$ and that $l_b$ is positive
semi\/{\rm -}\/definite. Then $l_a$ is positive semi\/{\rm -}\/definite.
\endproclaim

Indeed, both parameters of the signature $(s_+(\beta),s_-(\beta))$
of $l_\beta$ are lower semi-continuous functions. By the
conditions of the lemma, $s_+(\beta)+s_-(\beta)$ is constant on
$\left(a,b\right]$, and hence the same is true for
$s_{\pm}(\beta)$.

\demo{Proof  of Lemma {\rm \ref{reduce}}}
Since $\nrmint(\pi_t,0)$ is invertible, 
$\Inner(\pi_t,0)$ is a nondegenerate Hermitian form on
$I(\pi_t,0)$ for any $t$. Thus, the first statement follows from
Lemma \ref{elem}, after passing to any $K$-type.

To prove the second part, we will apply the discussion following
Proposition \ref{lqsetup} to the representations $\pi_t$. We may
identify all the induced spaces in (\ref{seq}) with the ones for
$t=0$ in the usual manner. The $K$-action will be independent of
$t$. We obtain a decomposition of the operator
$\ident(\pi_t,-\frac12)\nrmint(\pi_t,\frac12)$ defining the form
$\Inner(\pi_t,\frac12)$ as $C_t\circ B\circ A_t$ such that for
$t\ne0$ we have $\im(B\circ A_t)=\im(B)$ and $\ker(C_t\circ
B)=\ker(B)$. The crucial point is that the operator $B$ (denoted
by $R_{w_2}$ in (\ref{seq}) does not depend on $t$. Thus on each
$K$-type of $I(\pi_t,\frac12)$ the rank of $\Inner(\pi_t,\frac12)$
is equal to the rank of $B$, as long as $t\ne0$. Thus, we may
apply Lemma \ref{elem} to conclude the second statement of the lemma.
\enddemo

3.7. {\it The tempered Case}. We continue the proof of Lemma
\ref{claim2}. By virtue of the last section, we may assume that
$\pi$ is tempered. In this case, the representations $I(\pi,s)$
are irreducible for $0<s<\frac12$ by Lemma \ref{spn} and
Proposition \ref{char}. Thus $\Inner(\pi,s)$ is nondegenerate for
$0<s<\frac12$. We will show below that if $\Inner(\pi,\frac12)$ is
semi-definite then $\pi$ is of $G$-type. Then by Lemma
\ref{imply}, $\Inner(\pi,0)$ is definite. We may use Lemma
\ref{elem} on each $K$-type to conclude Lemma \ref{claim2}.

It remains to show that $\pi$ is of $G$-type if $\pi\in\Pi^\sdu$
is tempered and $\Inner(\pi,\frac12)$ is semi-definite. To shorten
notation, let $\pi_1=\pi^\orthpairs\times\pi^\symp\in\Pi^\sdu$ and
$\pi_2=\pi^\orth$ so that $\pi=\pi_1\times\pi_2$. Note that
$\pi_1$ is of $G$-type and hence $\herm(\pi_1,0)$ is a scalar by
Lemma \ref{imply}. Since $I(\pi_1,s)$ is irreducible for
$0<s<\frac12$ it follows from Lemma \ref{elem} that
\begin{equation} \label{definite}
\herm\left(\pi_1,\frac12\right)\hbox{ is semi-definite.}
\end{equation}
We need to show that $\pi_2=0$. Consider the family
$$
I\left(\pi_1\otimes\pi_2,\left(\frac12,\gamma\right)\right):=
I\left(\pi_1\nu^\frac12\times\pi_2\nu^\gamma,0\right).
$$
Let
$$
\herm'\left(\gamma\right):I\left(\pi_1\otimes\pi_2,\left(\frac12,\gamma\right)\right)
\rightarrow
I\left(\pi_1\otimes\pi_2,\left(-\frac12,-\gamma\right)\right)
$$
be the operator $\identt(\gamma)\circ\nrmintt(\gamma)$ where
$$
\nrmintt(\gamma):I\left(\pi_1\otimes\pi_2,\left(\frac12,\gamma\right)\right)
\longrightarrow
I\left(\opp{\pi_1}\otimes\opp{\pi_2},\left(-\frac12,-\gamma\right)\right)
$$
is the normalized intertwining operator
and
\begin{eqnarray*}
\identt(\gamma)=I\left(\ident_{\pi_1}\otimes\ident_{\pi_2},\left(-\frac12,-
\gamma\right)\right):
I\left(\opp{\pi_1}\otimes\opp{\pi_2},\left(-\frac12,-\gamma\right)\right)\\
\longrightarrow
I\left(\pi_1\otimes\pi_2,\left(-\frac12,-\gamma\right)\right).
\end{eqnarray*}
As usual we identify the underlying $K$-module of each family of
induced representations, so that it does not depend on $\gamma$.
The same argument as in Proposition \ref{lqsetup} with the
exponent $\gamma>0$ instead of $\frac12$ gives:
\specialnumber{3}\proclaim{Proposition} If $\gamma>0$ then
$I\left(\pi_1\otimes\pi_2,\left(\frac12,\gamma\right)\right)$
admits a Langlands quotient which is given by the image of
$\nrmintt(\gamma)$. It is isomorphic to
$\pi_2\nu^\gamma\times\LQ{\pi_1}$.
\endproclaim

The operator $\nrmintt(\gamma)$ can be written as the composition
of the intertwining maps
\begin{eqnarray*}
&&I\left(\pi_1\otimes\pi_2,\left(\frac12,\gamma\right)\right)\\
\longrightarrow&&I\left(\opp{\pi_1}\otimes\pi_2,\left(-\frac12,\gamma\right)
\right)\\
\xrightarrow{\id\sdp\nrmint(\pi_2,\gamma)}&&
 I\left(\opp{\pi_1}\otimes\opp{\pi_2},\left(-\frac12,-\gamma\right)\right)
\end{eqnarray*}
where the first map is
$$
I\left(\intgln\left(\pi_2\nu^\gamma,\opp{\pi_1}\nu^{-\frac12}\right),0\right
)\circ
\left(\id\sdp\nrmint\left(\pi_1,\frac12\right)\right)\circ
I\left(\intgln\left(\pi_1\nu^{\frac12},\pi_2\nu^\gamma\right),0\right).
$$
As before, the intermediate map $\id\sdp\nrmint(\pi_1,\frac12)$,
which does not depend on $\gamma$, already gives the Langlands
quotient as its image (on the full induced representation)
for  $\gamma>0$. Hence the rank of $\herm'(\gamma)$ on each
$K$-type is independent of~$\gamma$. Since $\herm'
(\frac12)=\herm(\frac12)$ we conclude by Lemma \ref{elem} that
$\herm'(0)$ is a semi-definite operator.

Now,
\begin{eqnarray*}
\identt(0)\circ\id\sdp\nrmint(\pi_2,0)&=&
\ident_{\pi_1}\nu^{-\frac12}\times\ident_{\pi_2}\sdp\triv\circ
\id\sdp\nrmint(\pi_2,0)\\&=&
\ident_{\pi_1}\nu^{-\frac12}\sdp\herm(\pi_2,0)=
\id\sdp\herm(\pi_2,0)\circ\ident_{\pi_1}\nu^{-\frac12}\sdp\id.
\end{eqnarray*}
Also, since $\pi_2\times\opp{\pi_1}\nu^{-\frac12}$ is irreducible,
$$
\left(\ident_{\pi_1}\nu^{-\frac12}\times\id\right)\circ\intgln\left(\pi_2,\opp{\pi_1}\nu^{-\frac12}\right)=
\intgln\left(\pi_2,\pi_1\nu^{-\frac12}\right)\circ\left(\id\times\ident_{\pi
_1}\nu^{-\frac12}\right)
$$
up to a scalar. All in all, $\herm'(0)$ is equal up to a scalar to
\begin{eqnarray*}
\id\sdp\herm\left(\pi_2,0\right)\circ\ident_{\pi_1}\nu^{-\frac12}\sdp\id\circ
I\left(\intgln\left(\pi_2,\opp{\pi_1}\nu^{-\frac12}\right),0\right)\circ
\left(\id\sdp\nrmint\left(\pi_1,\frac12\right)\right)\circ
M_1\\=\id\sdp\herm\left(\pi_2,0\right)\circ
I\left(\intgln\left(\pi_2,\pi_1\nu^{-\frac12}\right),0\right)\circ
\left(\id\sdp\herm\left(\pi_1,\frac12\right)\right)\circ M_1
\end{eqnarray*}
where $M_1=I(\intgln(\pi_1\nu^{\frac12},\pi_2),0)$. Note that by
the properties of the normalized intertwining operators
$I\left(\intgln\left(\pi_2,\pi_1\nu^{-\frac12}\right),0\right)$ is
the Hermitian dual of $M_1$ up to a scalar, and hence
$$
I\left(\intgln\left(\pi_2,\pi_1\nu^{-\frac12}\right),0\right)\circ
\left(\id\sdp\herm\left(\pi_1,\frac12\right)\right)\circ M_1
$$
is semi-definite. On the other hand, $\herm(\pi_2,0)$ is a
Hermitian involution. Thus, for $\herm'(0)$ to be semi-definite it
is necessary and sufficient that
\begin{equation} \label{contain}
M_1^{-1}\left(\ker\left(\id\sdp\,\nrmint\left(\pi_1,\frac12\right)\right)\right)\supset
\pi_1\nu^{\frac12}\sdp\,\Omega^\pm
\end{equation}
where $\Omega^\pm$ are the $\pm1$ eigenspaces of $\herm(\pi_2,0)$
on $I(\pi_2,0)$. Indeed, $\herm'(0)$ is semi-definite, of opposite
signs, on the subspaces $\pi_1\nu^\frac12\otimes\Omega^\pm$. We
will show that (\ref{contain}) is impossible if $\pi_2\ne0$. Let
$\omega$ be any irreducible constituent of $I(\pi_2,0)$. The
Langlands quotient of $\pi_1\nu^\frac12\sdp\,\omega$ is obtained
as the image of the corresponding intertwining operator (with
respect to a maximal parabolic of $G$)
$$
M_2:\pi_1\nu^\frac12\sdp\,\omega\rightarrow\opp{\pi_1}\nu^{-\frac12}\sdp\,\omega
$$
which is given by convergent integral. On the other hand, $M_2$ is
also the restriction to $\pi_1\nu^\frac12\sdp\omega$ of the
intertwining operator (with respect to a co-rank two parabolic
subgroup of $G$, but the same Weyl element)
$$
M_3:I\left(\pi_1\otimes\pi_2,\left(\frac12,0\right)\right)\rightarrow
I\left(\opp{\pi_1}\otimes\pi_2,\left(-\frac12,0\right)\right).
$$
Thus, we conclude that the image of $\pi_1\nu^\frac12\sdp\,\omega$
under $M_3$ is nonzero. On the other hand $M_3$ is obtained as
the composition of
\begin{eqnarray*}
 I\left(\pi_1\otimes\pi_2,\left(\frac12,0\right)\right)&\xrightarrow{M_1}&
I\left(\pi_2\otimes\pi_1,\left(0,\frac12\right)\right)\\
&   \xrightarrow{\id\sdp\nrmint\left(\pi_1,\frac12\right)}&
I\left(\pi_2\otimes\opp{\pi_1},\left(0,-\frac12\right)\right)\\
&\vxrightarrow
{I\left(\intgln\left(\pi_2,\opp{\pi_1}\nu^{-\frac12}\right),0\right)}&
I\left(\opp{\pi_1}\otimes\pi_2,\left(-\frac12,0\right)\right).
\end{eqnarray*}
Thus the left-hand side of (\ref{contain}) does not contain
$\pi_1\nu^\frac12\sdp\,\omega$ for any irreducible constituent of
$I(\pi_2,0)$. It remains to show that:

\specialnumber{11}\proclaim{Lemma}
$\Omega^\pm\ne0$ if $\pi_2\ne0$.
\endproclaim

\demo{Proof} This follow from the theory of $R$-groups (cf.\
\cite{Gol94}). Indeed, let
$\sigma=\rho_1\otimes\dots\otimes\rho_r$ considered as a
square-integrable representation of a Levi subgroup $L$ of $M={\rm
GL}_m$ and let $Q=LV$ be the corresponding standard parabolic
subgroup of $S_m$. Thus $\pi_2=\Ind_{Q\cap M}^M\sigma$. By our
conditions on $\sigma$, the $R$-group of $\sigma$ in $S_m$ is
isomorphic to
$$
W(\sigma)=\{w\in W/W_L:wLw^{-1}=L,w\sigma\isom\sigma\}.
$$
Thus any nontrivial element in $W(\sigma)$ gives rise to a
nonscalar intertwining operator $R_w$. Since the operator
$\herm(\pi_2,0)$ is up to a scalar $R_w$ for $w=w_0w_0^L$ we get
the result.
\enddemo

{\it Remark}.
Suppose that $\theta=\id$ and consider the following conditions on
a self-dual generic representation of ${\rm GL}_n$.
\begin{itemize}
\item[1.] $\pi$ is of $G$-type.
\item[2.] $\herm(\pi,0)$ is a scalar.
\item[3.] $I(\pi,\frac12)$ has a unitarizable quotient.
\item[4.] $\herm(\pi,\frac12)$ is semi-definite.
\end{itemize}
We remarked above that conditions 1 and 2
are equivalent. Similarly,   conditions 3 and
4 are equivalent. In the tempered case, all the
conditions are equivalent, although in general 3 is
stronger than 1. Any local component of a cuspidal
representation of $G$-type satisfies 3.
It seems that 3 reflects the fact that $\pi$
is a functorial image of a unitarizable representation of a
classical group.

\AuthorRefNames [CKPSS01]


\begin{references}

\bibitem{Art80}
\name{J.\  Arthur},
  A trace formula for reductive groups.\ {I}{I},  {A}pplications of a
  truncation operator,
   {\it Compositio Math\/}.\ {\bf 40} (1980), 87--121.

\bibitem{BG92}
\name{D.\ Bump} and \name{D.\  Ginzburg},
  Symmetric square ${L}$-functions on ${\rm {G}{L}}(r)$,
   {\it Ann.\ of Math\/}.\ {\bf 136} (1992), 137--205.

\bibitem{BH99}
\name{C.\ J.\ Bushnell} and \name{G.\  Henniart},
  Calculs de facteurs epsilon de paires pour ${\rm {G}{L}}\sb n$ sur un
  corps local.\ {I},
   {\it Bull.\ London Math.\ Soc\/}.\ {\bf 31} (1999), 534--542.

\bibitem{BM}
\name{E.\ M.\  Baruch} and \name{Z.\  Mao},
  Central critical value of automorphic $L$-functions,
  preprint.

\bibitem{BW00}
\name{A.~Borel} and \name{N.~Wallach},
   {\it Continuous Cohomology\/},  {\it Discrete Subgroups\/},
{\it and Representations
  of Reductive Groups}, 2nd edition,
  A.\ M.\ S., Providence, RI,  2000.

\bibitem{CI00}
\name{J.~B. Conrey} and \name{H.~Iwaniec},
  The cubic moment of central values of automorphic
${L}$-functions,
   {\it Ann.\ of Math\/}.\ {\bf 151} (2000), 1175--1216.

\bibitem{CKP-SS01}
\name{J.~W.\ Cogdell, H.~H.\ Kim, I.~I.\ Piatetski-Shapiro}, and
\name{F.~Shahidi},
  On lifting from classical groups to ${\rm {G}{L}}\sb {N}$,
   {\it Publ.\ Math.\ Inst.\ Hautes \'Etudes Sci\/}.\ {\bf 93}
(2001), 5--30.

\bibitem{CS98}
\name{W.\  Casselman} and \name{F.\ Shahidi},
  On irreducibility of standard modules for generic
representations.
   {\it Ann.\ Sci.\ \'Ecole Norm.\ Sup\/}.\  {\bf 31} (1998), 561--589.

\bibitem{Del76}
\name{P.\  Deligne},
  Les constantes locales de l'\'equation fonctionnelle de la fonction
  ${L}$ d'{A}rtin d'une repr\'esentation orthogonale,
   {\it Invent.\ Math\/}.\ {\bf 35} (1976), 299--316.

\bibitem{FQ73}
\name{A.~Fr{\"o}hlich} and \name{J.~Queyrut},
  On the functional equation of the {A}rtin ${L}$-function for
  characters of real representations,
   {\it Invent.\ Math\/}.\ {\bf 20} (1973), 125--138.

\bibitem{GJ72}
\name{R.\  Godement} and \name{H.\  Jacquet},
   {\it Zeta Functions of Simple Algebras},
  {\it Lecture Notes in Math\/}.\ {\bf 260},
  Springer-Verlag, New York, 1972.

\bibitem{GK75}
\name{I.~M.\ Gel$'$fand} and \name{D.~A.\ Kajdan},
  Representations of the group ${\rm {G}{L}}(n,{K})$ where ${K}$ is a
  local field,
  in  {\it Lie Groups and their Representations\/} (Proc.\
Summer School,
  Bolyai J\'anos Math.\ Soc\., Budapest, 1971),  95--118, Halsted, New
  York, 1975.

\bibitem{Gol94}
\name{D.\  Goldberg},
  Reducibility of induced representations for ${\rm {S}p}(2n)$ and
  ${\rm {S}{O}}(n)$,
   {\it Amer.\ J.\ Math\/}.\ {\bf 116} (1994), 1101--1151.

\bibitem{GP94}
\name{B.\ H.\ Gross} and \name{D.\ Prasad},
  On irreducible representations of ${\rm {S}{O}}\sb {2n+1}\times{\rm
  {S}{O}}\sb {2m}$,
   {\it Canad.\ J.\ Math\/}.\ {\bf 46} (1994), 930--950.

\bibitem{GRS01}
\name{D.\  Ginzburg, S.\  Rallis}, and \name{D.\ Soudry},
  Generic automorphic forms on\break ${\rm {S}{O}}(2n+1)$: functorial lift
to
  ${\rm {G}{L}}(2n)$, endoscopy, and base change,
   {\it Internat.\ Math.\ Res.\ Notices}  {\bf 14} (2001), 729--764.

\bibitem{Guo96}
\name{J.\  Guo},
  On the positivity of the central critical values of automorphic\break
  ${L}$-functions for ${\rm {G}{L}}(2)$,
   {\it Duke Math.\ J\/}.\ {\bf 83} (1996), 157--190.

\bibitem{IS00}
\name{H.\  Iwaniec} and \name{P.\  Sarnak},
  The nonvanishing of central values of automorphic\break ${L}$-functions
  and {L}andau-{S}iegel zeros,
   {\it Israel J.\ Math\/}.\ {\bf 120} (2000), 55--177.

\bibitem{Ivi01}
\name{A.\  Ivi{\'c}},
  On sums of {H}ecke series in short intervals,
   {\it J.\ Th{\rm \'{\it e}}or.\ Nombres Bordeaux} {\bf 13} (2001),
453--468.

\bibitem{Jan96}
\name{C.\  Jantzen},
  Reducibility of certain representations for symplectic and
  odd-orthogonal groups,
   {\it Compositio Math\/}.\ {\bf 104} (1996), 55--63.

\bibitem{JC01}
\name{H.\  Jacquet} and \name{N.\ Chen},
  Positivity of quadratic base change ${L}$-functions,
   {\it Bull.\ Soc.\ Math.\ France} {\bf 129} (2001), 33--90.

\bibitem{JP-SS83}
\name{H.~Jacquet, I.~I.\ Piatetskii-Shapiro}, and \name{J.~A.\ Shalika},
  Rankin-{S}elberg convolutions,
   {\it Amer.\ J.\ Math\/}.\ {\bf 105} (1983), 367--464.

\bibitem{JS}
\name{D.\  Jiang} and \name{D.\  Soudry},
  Generic representations and local Langlands reciprocity law for
  $p$-adic ${\rm SO}_{2n+1}$,
   {\it Contributions to Automorphic Forms, Geometry and Number Theory: Shalikafest
  2002},
  A Supplemental volume to the {\it Amer.\ J. Math.}, to appear, 2003.

\bibitem{JS81}
\name{H.~Jacquet} and \name{J.~A.\ Shalika},
  On {E}uler products and the classification of automorphic
  representations.\ {I},
   {\it Amer.\ J.\ Math\/}.\ {\bf 103} (1981), 499--558.

\bibitem{JS90a}
\bibline,
  Exterior square ${L}$-functions,
  in  {\it Automorphic Forms\/},  {\it Shimura Varieties\/},  {\it and
$L$-Functions\/},
   {\it Vol.\ II\/} (Ann Arbor, MI, 1988),  143--226. Academic Press,
  Boston, MA,
  1990.

\bibitem{JS90b}
\bibline,
  Rankin-{S}elberg convolutions: {A}rchimedean theory,
  in  {\it Festschrift in honor of I.\ I.\  Piatetski-Shapiro on the
  occasion of his sixtieth birthday\/},  {\it Part I\/} (Ramat Aviv, 1989),
  125--207, Weizmann, Jerusalem, 1990.

\bibitem{Kim02}
\name{H.\  Kim},
  On local $L$-functions and normalized intertwining operators,
   preprint, 2002.

\bibitem{KS88}
\name{C.~D.\  Keys} and \name{F.\  Shahidi},
  Artin ${L}$-functions and normalization of intertwining
operators,
   {\it Ann.\ Sci.\ \'Ecole Norm.\ Sup\/}.\ {\bf 21} (1988),
67--89.

\bibitem{KS93}
\name{S.\  Katok} and \name{P.\  Sarnak},
  Heegner points, cycles and {M}aass forms,
   {\it Israel J.\ Math\/}.\ {\bf 84} (1993), 193--227.

\bibitem{KZ81}
\name{W.~Kohnen} and \name{D.~Zagier},
  Values of ${L}$-series of modular forms at the center of the critical
  strip,
   {\it Invent.\ Math\/}.\ {\bf 64} (1981), 175--198.

\bibitem{Lan71}
\name{R.\ P.\ Langlands},
   {\it Euler Products},
  A James K. Whittemore Lecture in Mathematics given at Yale
  University, 1967, {\it Yale Mathematical Monographs\/} {\bf 1},
  Yale Univ.\ Press, New Haven, Conn., 1971.

\bibitem{Lan79}
\bibline,
  Automorphic representations, {S}himura varieties, and motives. {E}in
  {M}\"archen,
  in  {\it Automorphic Forms, Representations and
$L$-functions\/}, {\it Proc.\
  Sympos.\ Pure Math\/}.\ {\bf 33}  (Oregon State Univ., Corvallis, Ore.,
1977),
  205--246,\break A.\ M.\ S.,  Providence, R.I., 1979.

\bibitem{Lap02}
\name{E.\  Lapid},
  On the root number of representations of orthogonal type,
  {\it Compositio Math\/}., to appear.

\bibitem{Lap03}
\name{E.\  Lapid},
  On the nonnegativity of Rankin-Selberg ${L}$-functions at the center of
  symmetry,
   {\it Internat.\ Math.\ Res.\ Notices}  {\bf 2} (2003), 65--75.

\bibitem{Mui01}
\name{G.\  Mui{\'c}},
  A proof of {C}asselman-{S}hahidi's conjecture for quasi-split
  classical groups,
   {\it Canad.\ Math.\ Bull\/}.\ {\bf 44} (2001), 298--312.

\bibitem{MW89}
\name{C.~M{\oe{g}}lin} and \name{J.-L.\ Waldspurger},
  Le spectre r\'esiduel de ${\rm {G}{L}}(n)$,
   {\it Ann.\ Sci.\ {\rm \'{\it E}}cole Norm.\ Sup\/}.\ {\bf 22} (1989),
605--674.

\bibitem{MW95}
\bibline,
   {\it Spectral Decomposition and {E}isenstein Series},
  Cambridge Univ.\  Press, Cambridge, 1995,

\bibitem{PR99}
\name{D.\  Prasad} and \name{D.\  Ramakrishnan},
  On the global root numbers of ${\rm {G}{L}}(n)\times{\rm
{G}{L}}(m)$,
  in  {\it Automorphic Forms\/},  {\it Automorphic
Representations\/},  {\it and
  Arithmetic\/} (Fort Worth, TX, 1996), 311--330, A.\ M.\ S.,
  Providence, RI, 1999.

\bibitem{Ral87}
\name{S.\  Rallis},
   {\it ${L}$-Functions and the Oscillator Representation},
  Springer-Verlag, New York,  1987.

\bibitem{Sai95}
\name{T.\  Saito},
  The sign of the functional equation of the ${L}$-function of an
  orthogonal motive,
   {\it Invent.\ Math\/}.\ {\bf 120} (1995), 119--142.

\bibitem{Sha81}
\name{F.\  Shahidi},
  On certain ${L}$-functions,
   {\it Amer.\ J.\ Math\/}.\ {\bf 103} (1981), 297--355.

\bibitem{Sha84}
\bibline,
  Fourier transforms of intertwining operators and {P}lancherel
  measures for ${\rm {G}{L}}(n)$,
   {\it Amer.\ J.\ Math\/}.\ {\bf 106} (1984), 67--111.

\bibitem{Sha90a}
\bibline,
  On multiplicativity of local factors,
  in  {\it Festschrift in honor of\break I.\ I.\  Piatetski-Shapiro on the
  occasion of his sixtieth birthday\/},  {\it Part II\/}\break
  (Ramat Aviv, 1989),
  279--289, Weizmann, Jerusalem, 1990.

\bibitem{Sha90b}
\bibline,
  A proof of {L}anglands' conjecture on {P}lancherel measures;
  complementary series for  $p$-adic groups,
   {\it Ann.\ of Math\/}.\ {\bf 132} (1990), 273--330.

\bibitem{Sha91}
\bibline,
  Langlands' conjecture on {P}lancherel measures for $p$-adic
groups,
  in  {\it Harmonic Analysis on Reductive Groups\/} (Brunswick, ME,
1989),
   277--295. Birkh\"auser Boston, Boston, MA, 1991.

\bibitem{Sha92}
\bibline,
  Twisted endoscopy and reducibility of induced representations for\break
  $p$-adic groups,
   {\it Duke Math.\ J\/}.\ {\bf 66} (1992), 1--41.

\bibitem{Tad86}
\name{M.\  Tadi{\'c}},
  Classification of unitary representations in irreducible
  representations of general linear group (non{A}rchimedean case),
   {\it Ann.\ Sci.\ \'Ecole Norm.\ Sup\/}.\ {\bf 19} (1986), 335--382.

\bibitem{Tad98}
\bibline,
  On reducibility of parabolic induction,
   {\it Israel J.\ Math\/}.\ {\bf 107} (1998), 29--91.

\bibitem{Vog86}
\name{D.\ A.\  Vogan, Jr.},
  The unitary dual of ${\rm {G}{L}}(n)$ over an {A}rchimedean
field,
   {\it Invent.\ Math\/}.\ {\bf 83} (1986), 449--505.

\bibitem{Wal81}
J.-L.\ Waldspurger,
  Sur les coefficients de {F}ourier des formes modulaires de poids
  demi-entier,
   {\it J.\ Math.\ Pures Appl\/}.\ {\bf 60} (1981), 375--484.

\end{references}
\bye